\documentclass{article}

\usepackage{arxiv}
\usepackage[utf8]{inputenc}
\usepackage[T1]{fontenc}
\usepackage[numbers]{natbib}
\usepackage{amsmath,amssymb,amsfonts,amsthm,mathtools}
\usepackage{algorithm}
\usepackage{algpseudocode}
\usepackage{booktabs}
\usepackage{caption}
\usepackage{subcaption}
\usepackage{enumitem}
\usepackage{float}
\usepackage{graphicx}
\usepackage{microtype}
\usepackage{multicol}
\usepackage{needspace}
\usepackage{nicefrac}
\usepackage{titlesec}
\usepackage{xcolor}
\usepackage{url}
\usepackage{hyperref}
\usepackage{wrapfig}

\graphicspath{{./images/}}

\numberwithin{equation}{section}

\theoremstyle{definition}
\newtheorem{theorem}{Theorem}[section]

\newtheorem{lemma}[theorem]{Lemma}
\newtheorem{assumption}[theorem]{Assumption}
\newtheorem{remark}[theorem]{Remark}

\newcommand{\init}{\text{init}}
\newcommand{\new}{\text{new}}

\newcommand{\actor}{\text{actor}}
\newcommand{\critic}{\text{critic}}
\newcommand{\old}{\text{old}}

\title{Reinforcement-Learning-Guided Data-Driven Estimation of Spectral Properties of Stochastic Koopman Semigroups}

\author{Yuanchao Xu$^{1}$\thanks{Corresponding author email: xu.yuanchao.3a@kyoto-u.ac.jp}\enspace \enspace
Jing Liu$^{2}$\enspace
Weiping Ding$^{3}$\enspace
Zhongwei Shen$^{4}$\enspace
Isao Ishikawa$^{1}$\enspace\\
$^1$ Center for Science Adventure and Collaborative Research Advancement (SACRA)\\
Graduate School of Science,
Kyoto University\\
$^2$ Department of Biosystems Engineering, University of Manitoba\\
$^3$ School of Information Science and Technology, Nantong University\\
$^4$ Department of Mathematical and Statistical Science, University of Alberta\\
}

\begin{document}
\maketitle
\begin{abstract}
Koopman spectral analysis turns nonlinear stochastic dynamics into a linear evolution of observables and gives access to decay rates, oscillatory modes, and metastable behavior. In practice, however, EDMD, SDMD, and related estimators depend strongly on where the trajectory data are collected. If most trajectories start in regions that carry little spectral information, the leading eigenvalues and eigenfunctions can be poorly estimated even with a rich dictionary. We propose \emph{Reinforced SDMD}, a data-acquisition method that couples Stochastic Dynamic Mode Decomposition with reinforcement learning. The RL agent chooses trajectory-initialization regions, SDMD updates the Koopman approximation, and a spectral-consistency reward evaluates the estimated eigenpairs on the newly generated data. An exploration bonus is added to avoid repeatedly sampling only a small part of the state space. We test multi-armed bandits, DQN, and PPO on stochastic double-well, Duffing, and FitzHugh--Nagumo systems. The learned policies place more samples in regions that are useful for estimating the leading Koopman eigenpairs. We also give an error-propagation analysis showing how SDMD operator error enters the corresponding bandit, approximate value-iteration, and approximate policy-iteration bounds.
\end{abstract}

\section{Introduction}
Spectral analysis is a standard tool for studying long-time behavior, coherent structures, oscillations, and rare transitions in stochastic dynamical systems. The Koopman point of view studies the evolution of observables instead of states. Although the state dynamics may be nonlinear, the observables evolve linearly under the Koopman semigroup, and its eigenvalues and eigenfunctions describe intrinsic time scales and coherent features of the process~\cite{Koopman1931,koopman1932dynamical,mezic2005spectral}.

In applications, the Koopman spectrum must be estimated from finite trajectory data. Starting from Dynamic Mode Decomposition (DMD)~\cite{rowley2009spectral,tu2013dynamic}, methods such as EDMD~\cite{Williams_2015}, kernel EDMD~\cite{williams2015kernelkoopman}, Residual DMD~\cite{colbrook2024rigorous}, and SDMD for stochastic systems~\cite{xu2025datadrivenframeworkkoopmansemigroup} build finite-dimensional approximations of the Koopman operator, the Koopman semigroup, or the generator. Neural-network variants improve the choice of observables and help scale these methods to more complex data~\cite{pmlr-v119-azencot20a,baikonode,bai2025hierarchical,cheng2025machine,chretien2025using,gao2025coloke,li2017extended,lusch2018deep,yang2025tensor}.

A basic difficulty remains: the quality of a Koopman spectral estimate depends on the data distribution. Both hand-crafted dictionaries~\cite{ishikawa2024koopman,klus2020data,takeishi2017subspace,Williams_2015} and learned dictionaries~\cite{cheng2026informationshapeskoopmanrepresentation,1614066,li2017extended,doi:10.1137/18M1177846,xu2025reskoopnet} can fail to resolve important eigenfunctions if the sampled trajectories do not cover the relevant parts of state space. This issue is especially important for stochastic systems, where leading non-trivial eigenfunctions often describe slow mixing or rare transitions between metastable regions.

This paper treats trajectory generation as part of the learning problem. At each step, an agent chooses where to initialize new trajectories, SDMD updates the Koopman approximation, and the agent receives a reward based on the quality of the resulting spectral estimate. We call this method \emph{Reinforced SDMD}. The RL component is not a controller for the physical system. It is a data-acquisition policy whose actions select trajectory-initialization regions. We study three representative choices: multi-armed bandits~\cite{slivkins2019introduction}, Deep Q-Networks (DQN)~\cite{mnih2013playing}, and Proximal Policy Optimization (PPO)~\cite{schulman2017proximal}.

The experiments use a finite grid of initialization regions. This makes the bandit, DQN, and PPO results directly comparable and makes the learned sampling maps easy to interpret. The reward itself is not tied to a grid. A continuous-action version could let the policy output either an initial condition or the parameters of an initialization distribution, with the same SDMD estimator and spectral reward used for feedback.

The method is most useful when short trajectories can be generated or restarted from selected configurations, as in many simulations and some controlled experiments. Examples include computational fluid dynamics, molecular simulations, and neuronal models. If only a fixed passive dataset is available, then the same criterion can be used for adaptive subsampling or reweighting of trajectory segments, but not for closed-loop selection of new physical initial conditions.

The main contributions are as follows:
\begin{enumerate}
    \item We introduce \emph{Reinforced SDMD}, which combines SDMD with RL-based trajectory-initialization design for estimating spectra of stochastic Koopman semigroups under a limited data budget.
    \item We define a Koopman spectral-consistency reward and add a simple exploration bonus. The reward is tied to eigenpair consistency rather than short-horizon prediction error.
    \item We analyze how SDMD operator-estimation error propagates into bandit-, DQN-, and PPO-based data-acquisition schemes, and we test the method on double-well, stochastic Duffing, and stochastic FitzHugh--Nagumo systems.
\end{enumerate}

The paper is organized as follows. Section~2 reviews stochastic Koopman operators and SDMD. Section~3 presents Reinforced SDMD. Section~4 gives numerical results. Section~5 gives the error analysis. Section~6 concludes the paper.

\section{Preliminaries}
\subsection{Stochastic Koopman Operator}
The Koopman framework studies a dynamical system through observables. For stochastic dynamics, the evolution of an observable is given by conditional expectation, so the Koopman operators form a Markov semigroup.

Let $(X_t)_{t\ge0}$ be a continuous-time stochastic process on a probability space $(\Omega,\mathbb{P})$ satisfying
\begin{equation*}
    \mathrm{d}X_t=b(X_t)\,\mathrm{d}t+\sigma(X_t)\,\mathrm{d}W_t,
    \qquad X_0=x\in\mathcal{M},
\end{equation*}
where $\mathcal{M}\subseteq\mathbb{R}^d$, $b:\mathcal{M}\to\mathbb{R}^d$ is the drift, $\sigma:\mathcal{M}\to\mathbb{R}^{d\times r}$ is the diffusion matrix, and $W_t$ is an $r$-dimensional Wiener process.

Let $\rho$ be a probability measure on $\mathcal{M}$ and set $\mathcal{F}=L^2(\mathcal{M},\rho)$. For complex-valued observables, we use
\[
    \langle f,g\rangle_\rho=\int_{\mathcal{M}}\overline{f(x)}g(x)\,\mathrm{d}\rho(x).
\]
The stochastic Koopman semigroup $(\mathcal{K}^t)_{t\ge0}$ is defined by
\begin{equation*}
    (\mathcal{K}^t f)(x)=\mathbb{E}_{\mathbb{P}}[f(X_t)\mid X_0=x],
    \qquad f\in\mathcal{F}.
\end{equation*}
It gives the expected value of the observable after time $t$, starting from $x$.

The infinitesimal generator $\mathcal{A}$ is
\begin{equation*}
\mathcal{A}f=\lim_{t\downarrow0}\frac{\mathcal{K}^t f-f}{t},
\end{equation*}
with domain $\mathcal{D}(\mathcal{A})=\{f\in\mathcal{F}:\lim_{t\downarrow0}(\mathcal{K}^t f-f)/t\text{ exists in }\mathcal{F}\}$. For sufficiently smooth bounded functions, It\^{o}'s formula gives~\cite{pavliotis2016stochastic}
\begin{equation*}
    \mathcal{A}f
    =\sum_{i=1}^d b_i\frac{\partial f}{\partial x_i}
    +\frac12\sum_{i,j=1}^d(\sigma\sigma^\top)_{ij}
    \frac{\partial^2 f}{\partial x_i\partial x_j},
    \qquad f\in C_b^2(\mathcal{M}).
\end{equation*}

\begin{assumption}\label{ass:str_cont_semigroup_K}
Throughout the paper, $(\mathcal{K}^t)_{t\ge0}$ is a strongly continuous $C_0$-semigroup on $\mathcal{F}$.
\end{assumption}

Assumption~\ref{ass:str_cont_semigroup_K} implies that $\mathcal{D}(\mathcal{A})$ is dense in $\mathcal{F}$ and that $\mathcal{A}$ is closed~\cite{pazy2012semigroups}.

\begin{remark}\label{rem:spectral_mapping_eigenfunctions}
If $\phi\in\mathcal{D}(\mathcal{A})$ satisfies $\mathcal{A}\phi=\lambda\phi$, then
\begin{equation*}
    \mathcal{K}^t\phi=e^{t\lambda}\phi.
\end{equation*}
Thus the corresponding semigroup eigenvalue is $\mu=e^{t\lambda}$. The converse direction requires the usual spectral-mapping assumptions and is not used by the algorithms below.
\end{remark}

\subsection{Stochastic Dynamic Mode Decomposition (SDMD)}
EDMD approximates the Koopman operator on a finite-dimensional dictionary space. SDMD~\cite{xu2025datadrivenframeworkkoopmansemigroup} adapts this idea to stochastic differential equations by using the generator and a first-order stochastic Taylor expansion of the semigroup.

Let $\{\psi_1,\ldots,\psi_N\}\subset\mathcal{D}(\mathcal{A})$ and $\mathcal{F}_N=\operatorname{span}\{\psi_1,\ldots,\psi_N\}$. Given i.i.d. points $\{x_k\}_{k=1}^m$ sampled from $\rho$, define
\begin{equation*}
    [\Psi_X]_{kj}=\psi_j(x_k),
    \qquad
    [\Psi'_X]_{kj}=(\mathcal{A}\psi_j)(x_k),
\end{equation*}
and
\begin{equation*}
    \widehat G=\frac1m\Psi_X^*\Psi_X,
    \qquad
    \widehat H=\frac1m\Psi_X^*\Psi'_X.
\end{equation*}
Here $*$ denotes conjugate transpose. In computations we use $(\widehat G+\gamma_{\mathrm{reg}}I)^{-1}$ when $\widehat G$ is ill-conditioned. The first-order SDMD approximation of the lag-$\Delta t$ Koopman semigroup is
\begin{equation*}
    \widehat K_{N,\Delta t,m}=I+\Delta t\,\widehat G^{-1}\widehat H.
\end{equation*}
The ordered convergence limits are $m\to\infty$, $\Delta t\to0$, and then $N\to\infty$; see~\cite[Section~4]{xu2025datadrivenframeworkkoopmansemigroup}.

When the dictionary is learned, we parameterize it by a neural network $\Psi(x;\theta)$ and use a Koopman-consistency loss such as
\begin{equation*}
    J(\theta)
    =\|\Psi_Y(\theta)-\Psi_X(\theta)\widehat K_{N,\Delta t,m}(\theta)\|_F^2
    +\lambda_{\mathrm{reg}}\mathcal{R}_{\widehat K},
\end{equation*}
where $\Psi_Y$ is the dictionary evaluated on lagged data, $\mathcal{R}_{\widehat K}$ is a regularizer, and $\lambda_{\mathrm{reg}}>0$.

\section{Methodology: Reinforced Koopman Spectral Learning}
\begin{figure*}[!t]
    \centering
    \includegraphics[width=0.99\linewidth]{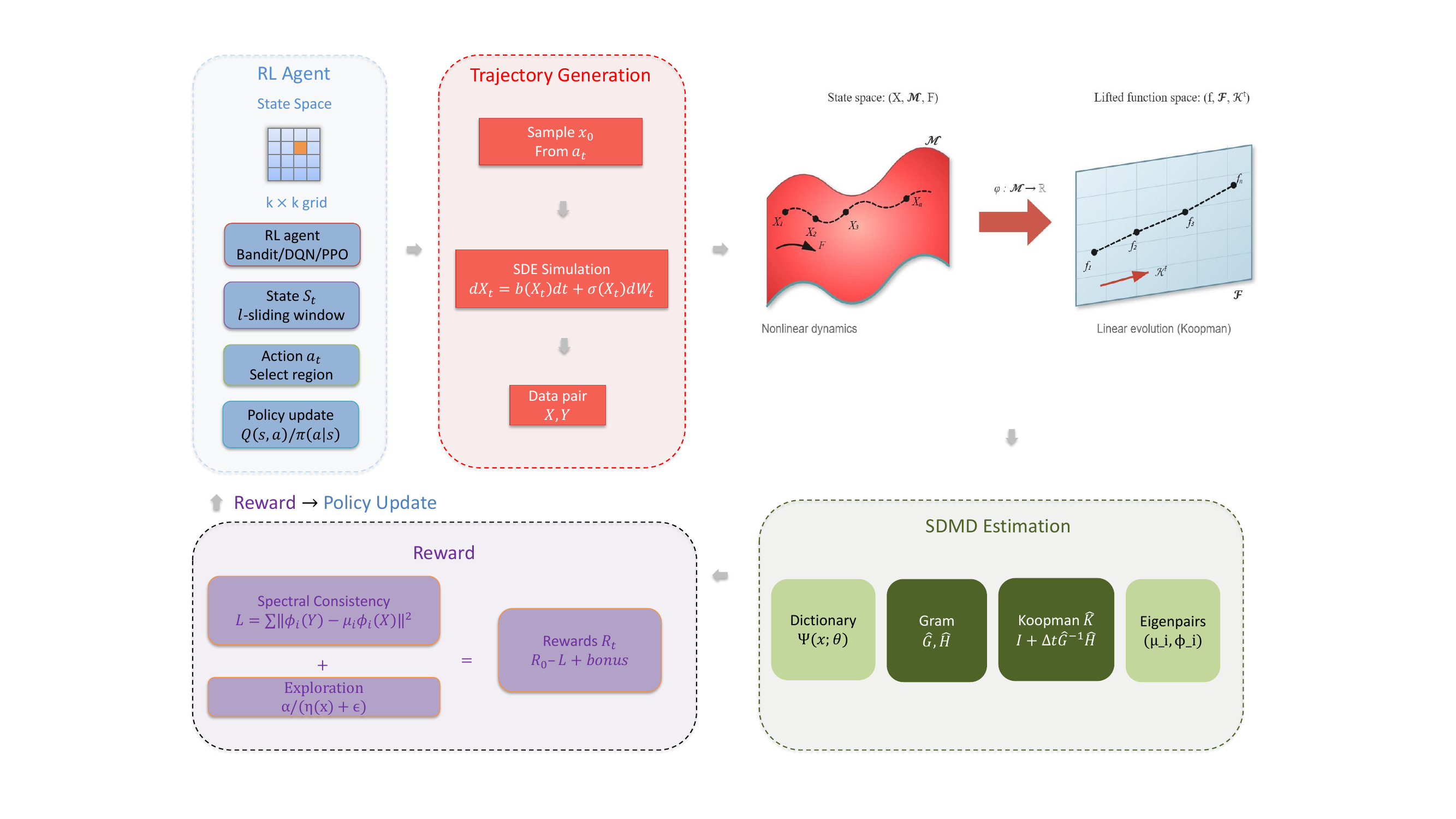}
    \caption{Workflow of Reinforced SDMD. The RL agent selects trajectory-initialization regions, SDMD updates the Koopman approximation, and a spectral-consistency reward guides future data acquisition.}
    \label{fig:flowchart}
\end{figure*}

Reinforced SDMD casts data collection as a sequential experiment-design problem. At iteration $t$, an agent selects a region of state space, one or more short trajectories are initialized in that region, SDMD updates the Koopman approximation, and the agent receives a reward based on spectral consistency. Figure~\ref{fig:flowchart} summarizes the workflow.

\subsection{Review of RL Algorithms}
\begin{figure}
    \centering
    \includegraphics[width=0.75\linewidth]{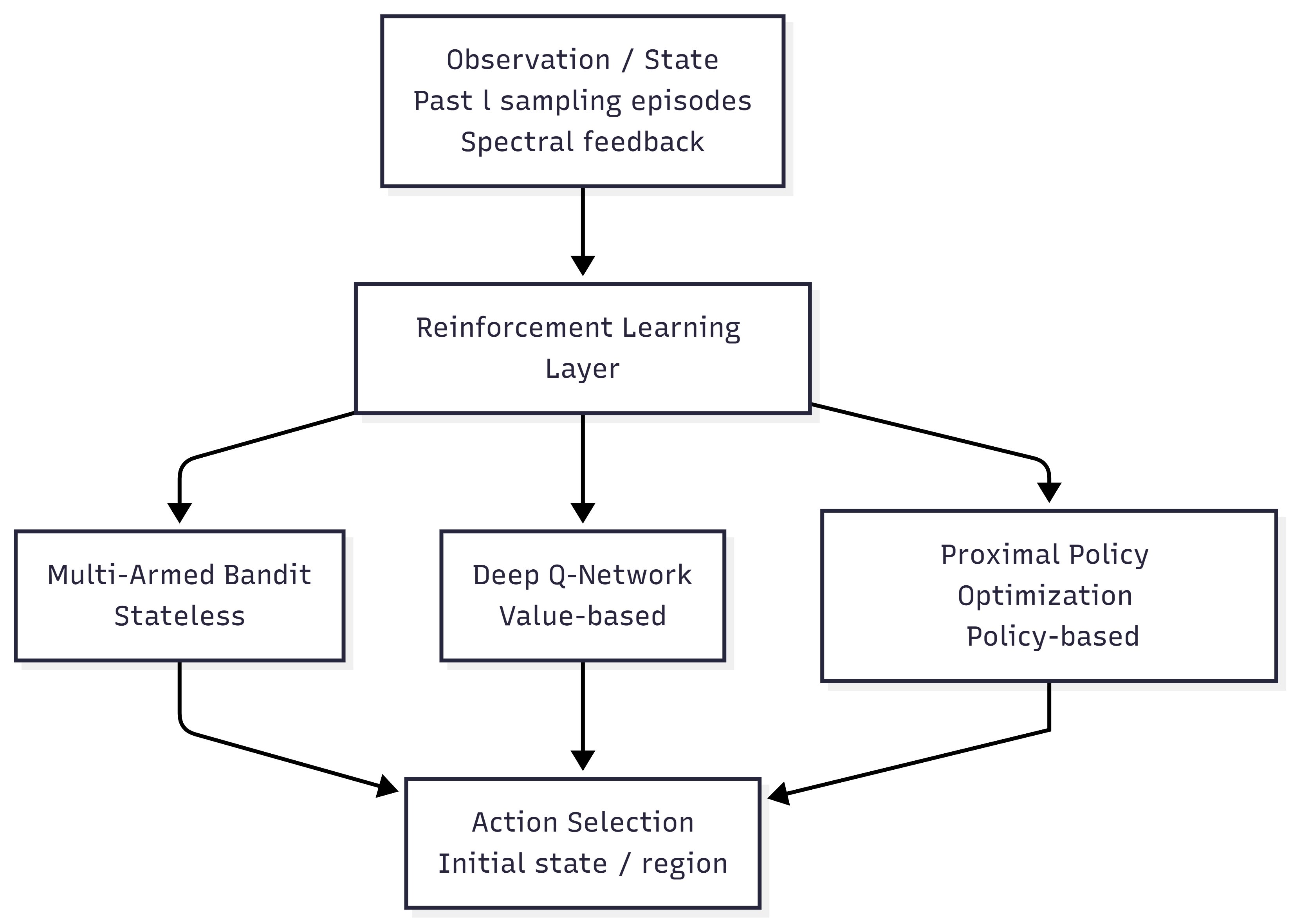}
    \caption{Three RL instantiations of Reinforced SDMD: an $\varepsilon$-greedy bandit, DQN, and PPO.}
    \label{fig:flowchart2}
\end{figure}

We use three standard RL schemes. Figure~\ref{fig:flowchart2} shows how they fit into the same Reinforced SDMD loop.

\textbf{Multi-armed bandit with $\varepsilon$-greedy policy.}
In the bandit setting~\cite{slivkins2019introduction}, each action selects one spatial region for trajectory initialization. The agent maintains an empirical value $Q(a)$ for each region and uses an $\varepsilon$-greedy rule to balance exploitation and random exploration.

\textbf{Deep Q-Network (DQN).}
DQN~\cite{mnih2013playing} uses a state variable. Here the state is a short history of recently selected initialization points. A neural network approximates $Q(S_t,a_t;\theta)$ and is trained by temporal-difference updates with experience replay and a target network.

\textbf{Proximal Policy Optimization (PPO).}
PPO~\cite{schulman2017proximal} is an actor--critic method. The actor learns a stochastic policy $\pi_\theta(a\mid s)$ over initialization regions, and the critic estimates the value function. The clipped PPO objective helps keep policy updates stable when the SDMD reward is noisy.

\begin{remark}
    Appendix~\ref{app:rl_review} gives additional algorithmic details.
\end{remark}

\subsection{Reinforced SDMD}
The connection between RL and SDMD is defined by the action space, the state representation, and the reward.

\textbf{Action space.}
An action $a_t$ selects a subset of the state space, and the next initialization point $x_{\new}$ is sampled from that subset. In the experiments, the subsets are grid cells. This finite action space makes the three RL methods easy to compare. The same idea can be used with continuous actions: the policy may output an initial point in $\mathcal{M}$ or the parameters of a distribution over initial points. In that case, only the policy class and the exploration term need to change.

\textbf{State representation.}
For DQN and PPO, the RL state is a sliding window
\[
    S_t=(x_{t-\ell+1},\ldots,x_t)\in\mathcal{M}^{\ell}
\]
of the last $\ell$ initialization points. This gives the agent a simple memory of where it has recently sampled.

\textbf{Reward formulation.}
Let $(X,Y)=\{(x_j,y_j)\}_{j=1}^{m}$ be paired snapshots with lag time $\Delta t$, where $y_j$ is one stochastic continuation from $x_j$. Suppose SDMD returns approximate eigenpairs $(\mu_i,\phi_i)$. For an exact stochastic Koopman eigenfunction,
\begin{equation*}
    (\mathcal K^{\Delta t}\phi_i)(x)
    =\mathbb{E}[\phi_i(Y)\mid X=x]
    =\mu_i\phi_i(x).
\end{equation*}
The identity is therefore a conditional-mean relation, not a pointwise equality for each sampled transition. In the experiments we use the empirical residual
\begin{equation*}
    \mathcal{L}_{\mathrm{SC}}(X,Y)
    \coloneqq
    \sum_{i=1}^{N_e}\frac{1}{m}\sum_{j=1}^{m}
    \left|\phi_i(y_j)-\mu_i\phi_i(x_j)\right|^2,
\end{equation*}
where $N_e$ is the number of eigenpairs included in the reward. This residual is a noisy proxy for spectral consistency. Indeed,
\begin{equation*}
\mathbb{E}\!\left[\left|\phi_i(Y)-\mu_i\phi_i(x)\right|^2\mid X=x\right]
=\left|(\mathcal K^{\Delta t}\phi_i)(x)-\mu_i\phi_i(x)\right|^2
+\operatorname{Var}\!\left(\phi_i(Y)\mid X=x\right),
\end{equation*}
where for complex-valued observables we use $\operatorname{Var}(Z)=\mathbb{E}|Z-\mathbb{E}Z|^2$. If repeated continuations from the same point or local averaging within a region are available, $\phi_i(y_j)$ can be replaced by a Monte Carlo estimate of $\mathbb{E}[\phi_i(Y)\mid X=x_j]$, reducing the variance term.

To avoid sampling only one region, we add a density-based exploration bonus. Let $\eta(x_{\new})$ be a kernel-density estimate computed from previous initialization points. The reward is
\begin{equation}\label{eq:reward}
    R_t
    \coloneqq
    R_0-\mathcal{L}_{\mathrm{SC}}(X,Y)
    +\alpha_{\mathrm{exp}}\frac{1}{\eta(x_{\new})+\varepsilon_{\mathrm{dens}}},
\end{equation}
where $R_0$ is a baseline, $\alpha_{\mathrm{exp}}>0$ controls exploration, and $\varepsilon_{\mathrm{dens}}>0$ prevents division by zero.

\textbf{Bandit-SDMD.}
The bandit version treats each grid cell as an arm. The learned $Q$-values give a reward map over initialization regions.
\begin{algorithm}
\caption{Multi-Armed Bandit with SDMD}
\label{alg:bandit_sdmd}
\begin{algorithmic}[1]
\Require Grid size $k$, number of steps $T_{\max}$, exploration schedule $\varepsilon_t$, initial value $Q_{\init}$
\State Initialize $Q(a)\leftarrow Q_{\init}$ and $N(a)\leftarrow0$ for all $a\in\{0,1,\ldots,k^2-1\}$.
\For{$t=1$ to $T_{\max}$}
    \State Select $a_t$ using the $\varepsilon_t$-greedy rule.
    \State Sample $x_{\new}$ from the grid cell indexed by $a_t$ and generate paired data $(X,Y)$.
    \State Compute SDMD eigenpairs and reward $R_t$ by Eq.~\eqref{eq:reward}.
    \State Update $N(a_t)\leftarrow N(a_t)+1$.
    \State Update $Q(a_t)\leftarrow Q(a_t)+N(a_t)^{-1}(R_t-Q(a_t))$.
\EndFor
\end{algorithmic}
\end{algorithm}

\textbf{DQN-SDMD.}
DQN-SDMD learns a state-dependent action-value function. After an action is selected, the new trajectory is combined with the trajectories associated with the current sliding-window state. SDMD computes the reward, and the DQN is trained using the temporal-difference target
\begin{equation*}
    R_t+\gamma\max_{a'}Q(S_{t+1},a';\theta^-).
\end{equation*}
The loss is the Huber loss between this target and $Q(S_t,a_t;\theta)$. The pseudocode is given in Algorithm~\ref{alg:dqn_sdmd}; implementation details are in Appendix~\ref{app:dqn}.
\begin{algorithm}
\caption{DQN with SDMD}
\label{alg:dqn_sdmd}
\begin{algorithmic}[1]
\Require State length $\ell$, replay capacity $N_{\mathrm{replay}}$, discount $\gamma$, soft-update parameter $\tau$, exploration schedule $\varepsilon_t$
\State Initialize policy network $Q(s,a;\theta)$ and target network $Q(s,a;\theta^-)$ with $\theta^-\leftarrow\theta$.
\State Initialize replay memory $\mathcal{D}$ and initial state $S_0$ using $\ell$ random initialization points.
\For{$t=1$ to $T_{\max}$}
    \State Select $a_t$ from $S_t$ using an $\varepsilon_t$-greedy rule.
    \State Generate new data from the selected region, compute SDMD eigenpairs, and compute $R_t$ by Eq.~\eqref{eq:reward}.
    \State Construct $S_{t+1}$ by the sliding-window update and store $(S_t,a_t,R_t,S_{t+1})$ in $\mathcal{D}$.
    \State Sample a minibatch from $\mathcal{D}$ and set $y_j=R_j+\gamma\max_{a'}Q(S_{j+1},a';\theta^-)$.
    \State Update $\theta$ by minimizing the Huber loss between $y_j$ and $Q(S_j,a_j;\theta)$.
    \State Soft-update $\theta^-\leftarrow\tau\theta+(1-\tau)\theta^-$.
\EndFor
\end{algorithmic}
\end{algorithm}

\textbf{PPO-SDMD.}
PPO-SDMD uses the same reward but learns a stochastic actor directly. The clipped surrogate objective is
\begin{equation}\label{eq:clip_actor}
    \mathcal{L}^{\mathrm{CLIP}}_{\mathrm{PPO}}(\theta_{\actor})
    \coloneqq
    \widehat{\mathbb{E}}_t
    \left[
    \min\left(
    r_t(\theta_{\actor})\widehat{A}_t,
    \operatorname{clip}\bigl(r_t(\theta_{\actor}),1-\varepsilon_{\mathrm{clip}},1+\varepsilon_{\mathrm{clip}}\bigr)\widehat{A}_t
    \right)
    \right],
\end{equation}
where $r_t$ is the policy likelihood ratio and $\widehat A_t$ is an advantage estimate.
\begin{algorithm}
\caption{PPO with SDMD}
\label{alg:ppo_sdmd}
\begin{algorithmic}[1]
\Require State length $\ell$, batch size $N_{\mathrm{batch}}$, discount $\gamma$, GAE parameter $\lambda$, clipping parameter $\varepsilon_{\mathrm{clip}}$, epochs $K$
\State Initialize actor $\pi(a\mid s;\theta_{\actor})$ and critic $V(s;\theta_{\critic})$.
\For{iteration $=1,\ldots,M$}
    \State Collect $N_{\mathrm{batch}}$ transitions using the current policy $\pi(\cdot\mid S_t;\theta_{\actor,\old})$.
    \State For each transition, generate data from the selected region, compute the SDMD reward, and update the sliding-window state.
    \State Compute advantage estimates $\widehat A_t$ by GAE.
    \For{epoch $=1,\ldots,K$}
        \State Update the actor by maximizing Eq.~\eqref{eq:clip_actor}.
        \State Update the critic by minimizing the value-function loss.
    \EndFor
    \State Set $\theta_{\actor,\old}\leftarrow\theta_{\actor}$.
\EndFor
\end{algorithmic}
\end{algorithm}

\begin{remark}
Several related approaches use RL, optimal experimental design, active learning, or adaptive sampling for scientific data acquisition. Here the reward is specifically tied to Koopman spectral consistency, rather than to prediction error or a task-specific control return~\cite{fannjiang2020autofocused,grayeli2024symbolic,hsu2010algorithms,shen2025variational,zhao2024policy,zhang2024adaptive}.
\end{remark}

\section{Numerical Experiments}\label{sec:experiments}
We test Reinforced SDMD on three stochastic systems: a double-well potential, a stochastic Duffing oscillator, and a stochastic FitzHugh--Nagumo model. These examples cover metastability, bistability with noise-induced switching, and fast--slow relaxation dynamics. Unless otherwise stated, the experiments were run on a server with an NVIDIA RTX 6000 Ada Generation GPU (48 GB VRAM), 16 vCPUs, and 62 GB RAM; each reported configuration finished within 24 hours.

\subsection{Double Well Potential System}

\begin{figure}[!htb]
    \centering

    \begin{subfigure}{0.33\textwidth}
        \centering
        \includegraphics[width=\linewidth]{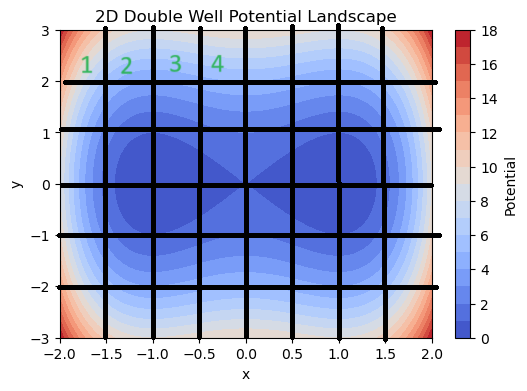}
        \caption{2D double-well potential landscape}
        \label{fig:2d_potential}
    \end{subfigure}
    \hspace{0.2cm}
    \begin{subfigure}{0.39\textwidth}
        \centering
        \includegraphics[width=\linewidth]{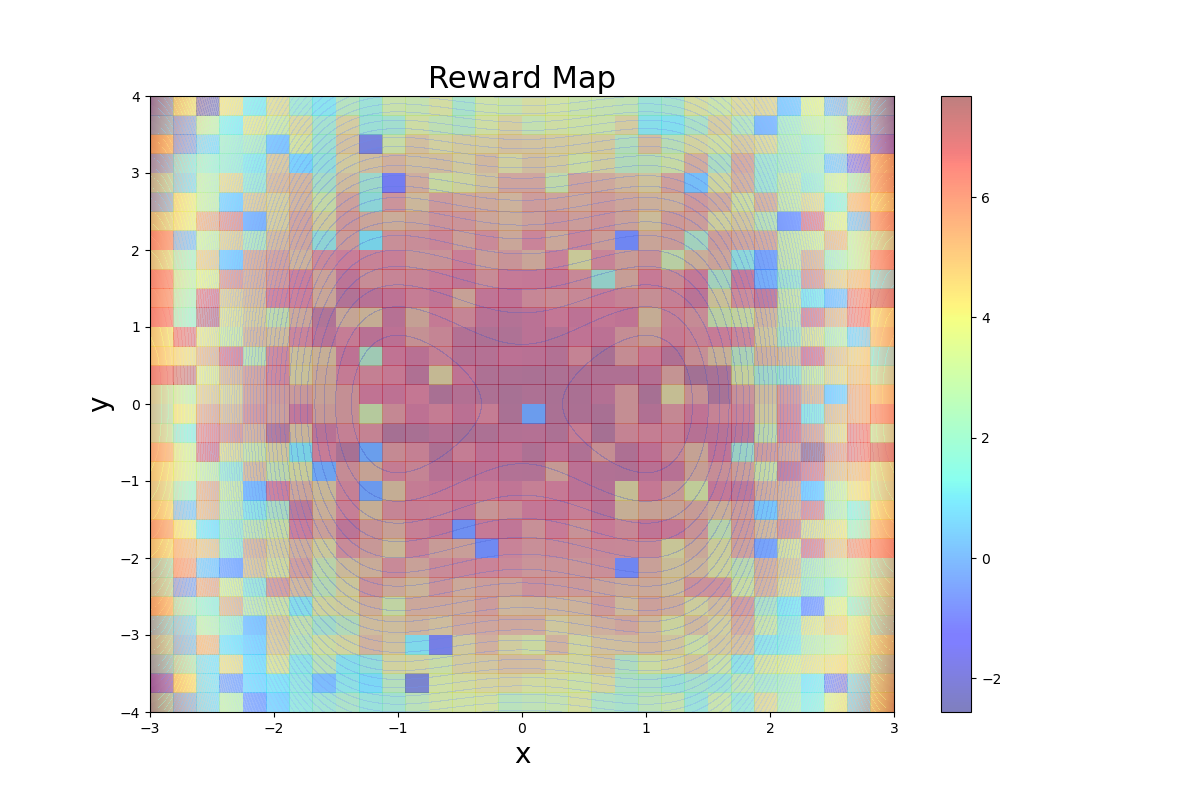}
        \caption{Learned reward map}
        \label{fig:reward_map}
    \end{subfigure}
    \caption{Double-well potential landscape and the learned Bandit-SDMD reward map.}
    \label{fig:comb_potential_reward}
\end{figure}

We first consider the two-dimensional double-well potential system
\begin{equation*}
    \mathrm{d}X_t=-\nabla V(X_t)\,\mathrm{d}t+\Sigma\,\mathrm{d}W_t,
\end{equation*}
where $X_t=(x_t,y_t)$, $V(x,y)=(x^2-1)^2+y^2$, and $\Sigma=1.09 I_2$. The domain $[-3,3]\times[-4,4]$ is discretized into a $32\times32$ grid, producing $1024$ possible actions. Each action selects one grid cell from which a new trajectory is initialized, as illustrated in Figure~\ref{fig:2d_potential}.

The reward follows Eq.~\eqref{eq:reward} and combines spectral consistency with an exploration bonus. We use $\varepsilon=0.35$ for the $\varepsilon$-greedy policy and $\alpha_{\mathrm{exp}}=0.15$ for the exploration coefficient. Each selected action generates a 1000-step trajectory, which is processed by SDMD with a neural-network dictionary; the bandit $Q$-values are then updated by sample averaging.

\begin{figure}[!htb]
    \centering
    \begin{subfigure}{0.45\linewidth}
        \centering
        \includegraphics[width=\linewidth]{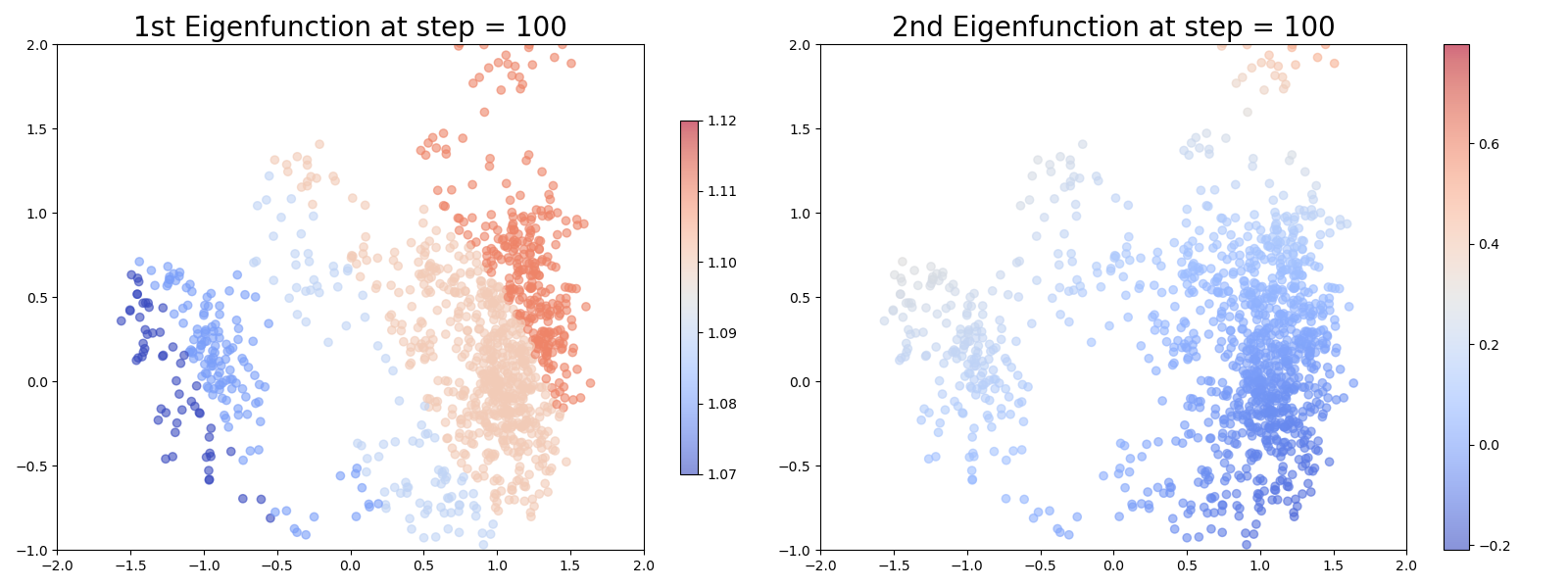}
        \caption{Step 100}
    \end{subfigure}\hspace{0.0cm}
    \begin{subfigure}{0.45\linewidth}
        \centering
        \includegraphics[width=\linewidth]{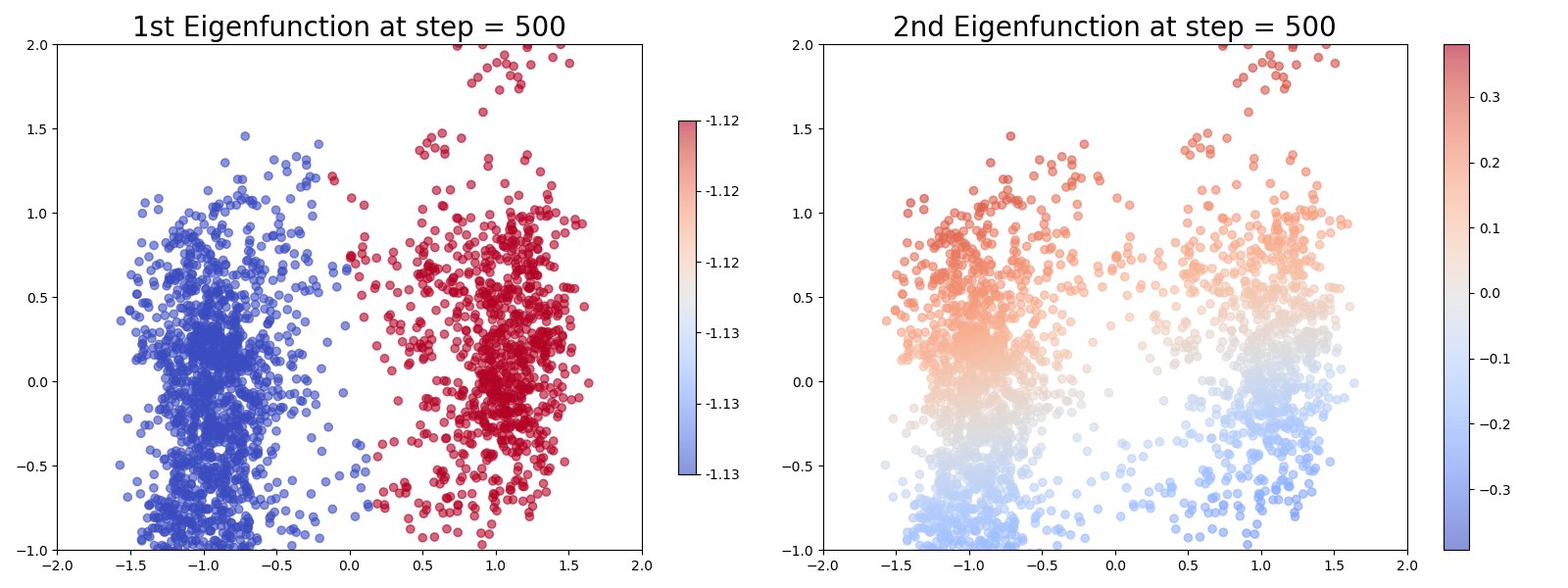}
        \caption{Step 500}
    \end{subfigure}

    \vspace{0.1cm}

    \begin{subfigure}{0.45\linewidth}
        \centering
        \includegraphics[width=\linewidth]{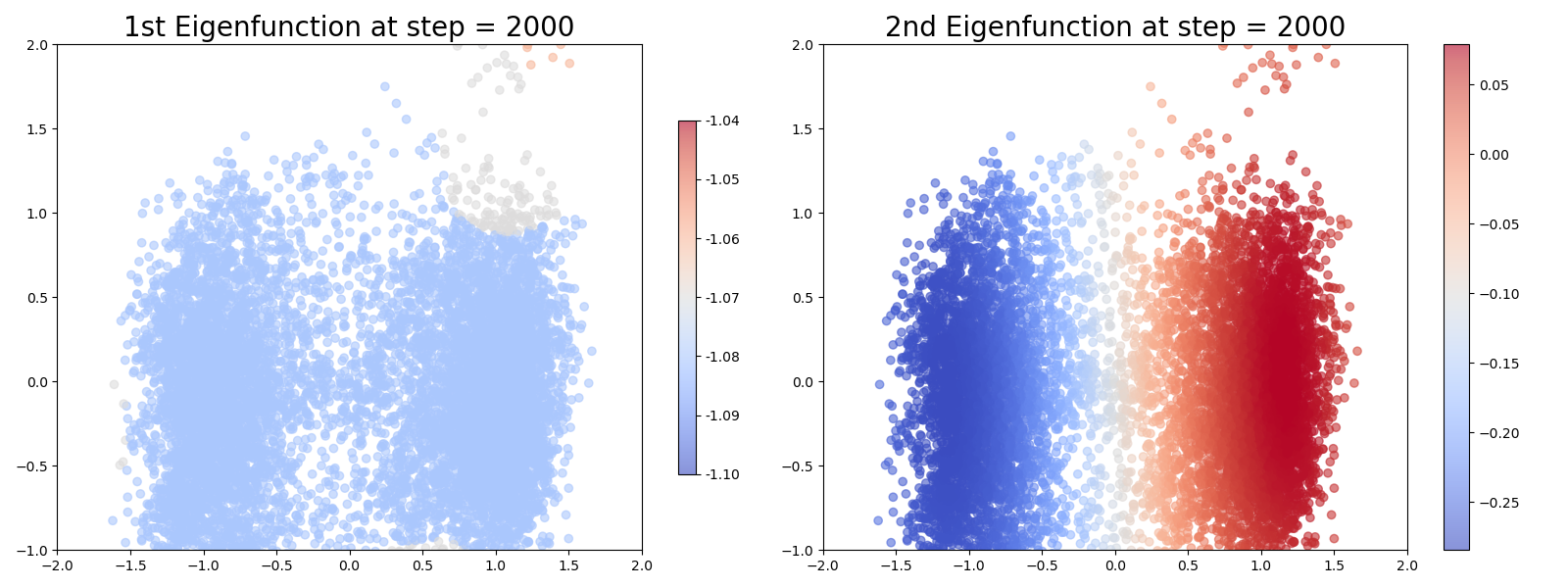}
        \caption{Step 2000}
    \end{subfigure}\hspace{0.0cm}
    \begin{subfigure}{0.45\linewidth}
        \centering
        \includegraphics[width=\linewidth]{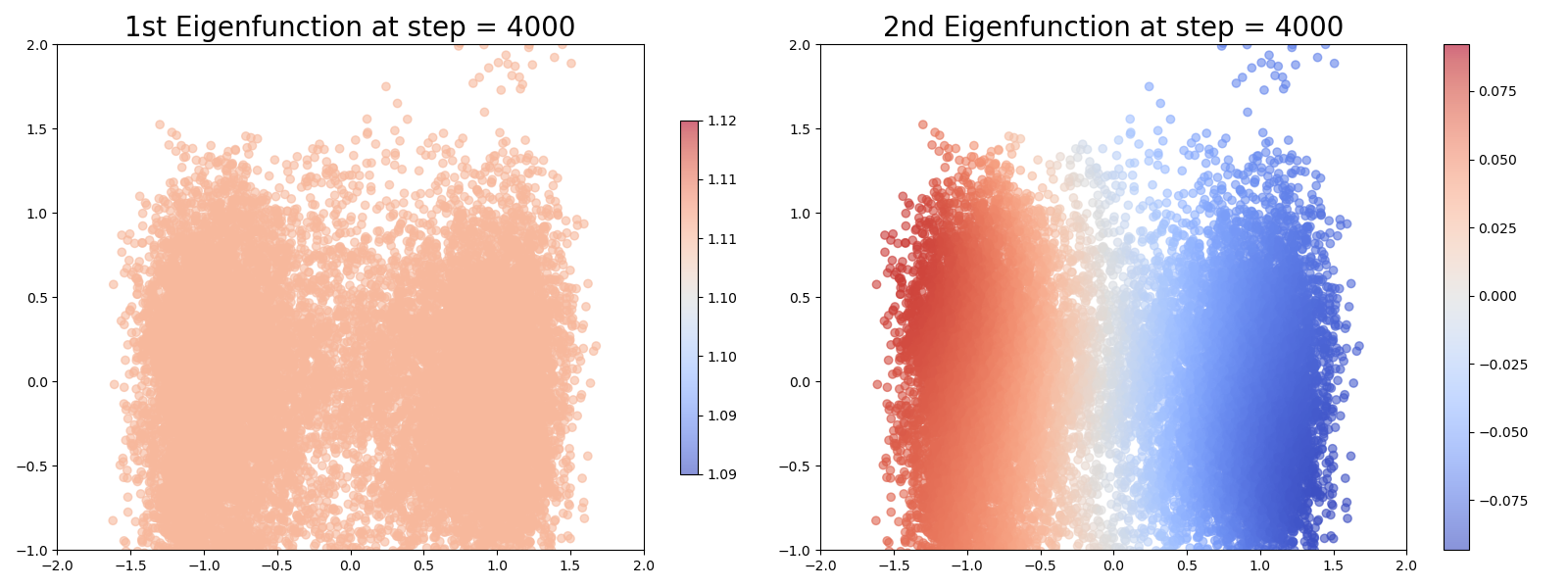}
        \caption{Step 4000}
    \end{subfigure}

    \caption{First and second Koopman eigenfunctions of the double-well potential system at different training steps.}
    \label{fig:double_well_bandit_steps}
\end{figure}

After 4000 steps, the learned reward map places more weight on the central part of the domain that contains the two wells and the transition region; see Figure~\ref{fig:reward_map}. This is the part of the state space most relevant for the leading metastable eigenfunctions. Figure~\ref{fig:double_well_bandit_steps} shows the same behavior at the level of eigenfunctions: the first eigenfunction approaches the invariant mode, while the second separates the two wells by taking different signs in the two basins.

\subsection{Stochastic Duffing Oscillator}

Consider the stochastic damped Duffing oscillator
\begin{equation*}
    \begin{cases}
    \mathrm d x(t) = v(t)\,\mathrm d t,\\
    \mathrm d v(t) = \bigl(-\delta v(t)-\alpha x(t)-\beta x(t)^{3}\bigr)\,\mathrm d t
                +\sigma\,\mathrm d W(t),
    \end{cases}
\end{equation*}




\begin{wrapfigure}{r}{0.4\linewidth}
    \centering
    \includegraphics[width=0.85\linewidth]{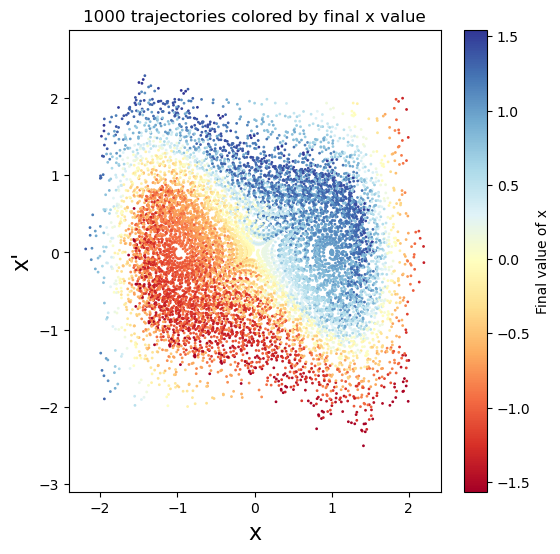}
    \caption{Stochastic Duffing oscillator: phase-space trajectory.}
    \label{fig:duffing_phase_portrait}
    \vspace{-3em}
\end{wrapfigure}
\begin{figure}[!htb]

    \centering
    \def\figscale{0.7}

    \begin{subfigure}{\figscale\linewidth}
        \centering
        \includegraphics[width=\linewidth]{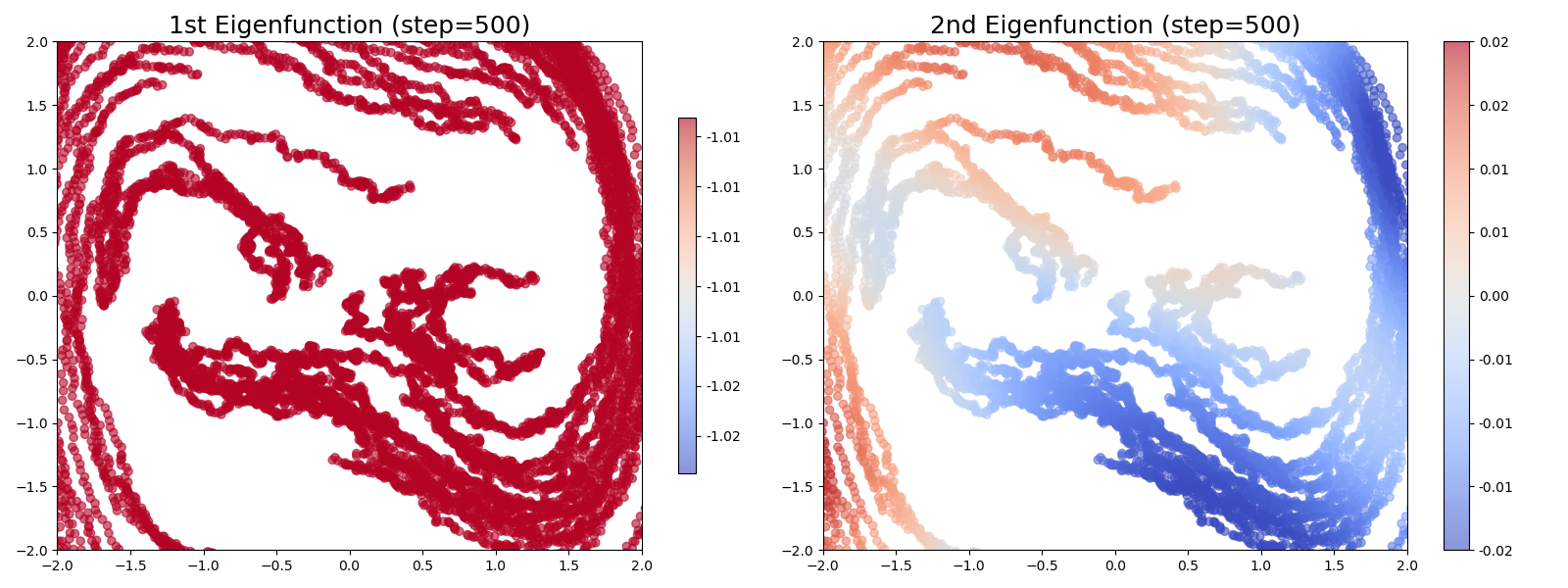}
        \caption{Step 500}
        \label{fig:duffing_dqn_step_500}
    \end{subfigure}

    \vspace{0.5em}

    \begin{subfigure}{\figscale\linewidth}
        \centering
        \includegraphics[width=\linewidth]{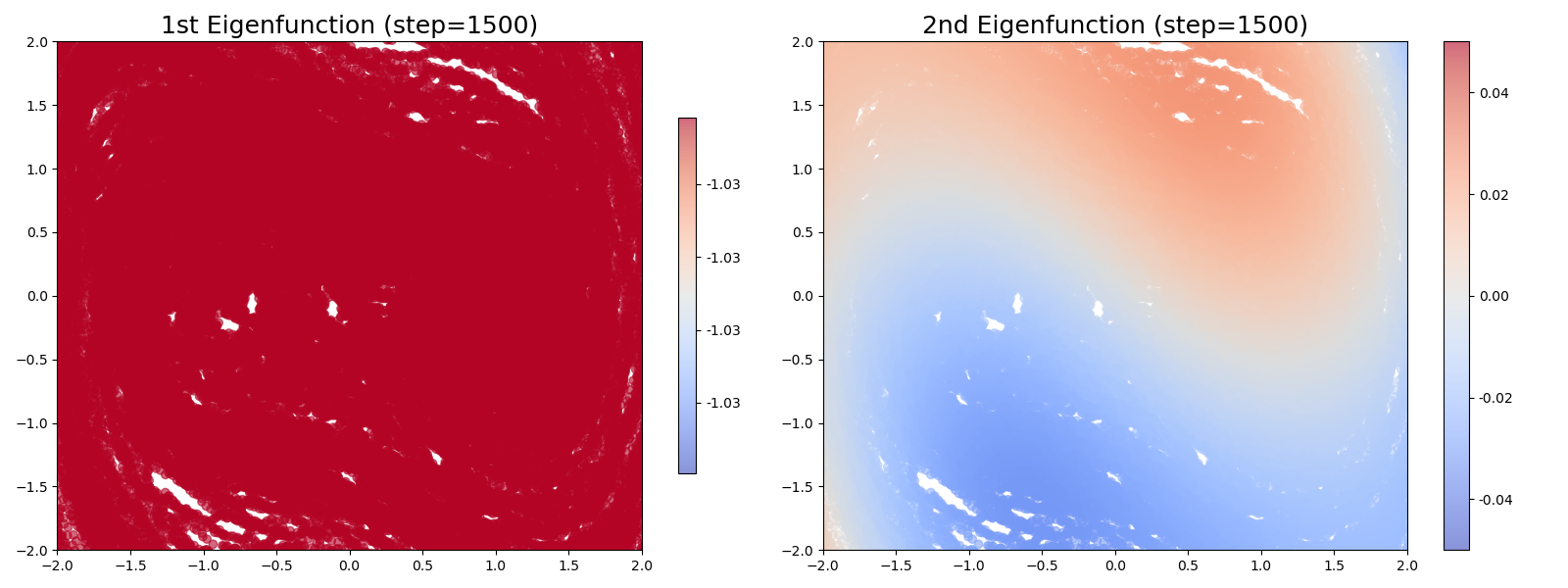}
        \caption{Step 1500}
        \label{fig:duffing_dqn_step_1500}
    \end{subfigure}

    \vspace{0.5em}

    \begin{subfigure}{\figscale\linewidth}
        \centering
        \includegraphics[width=\linewidth]{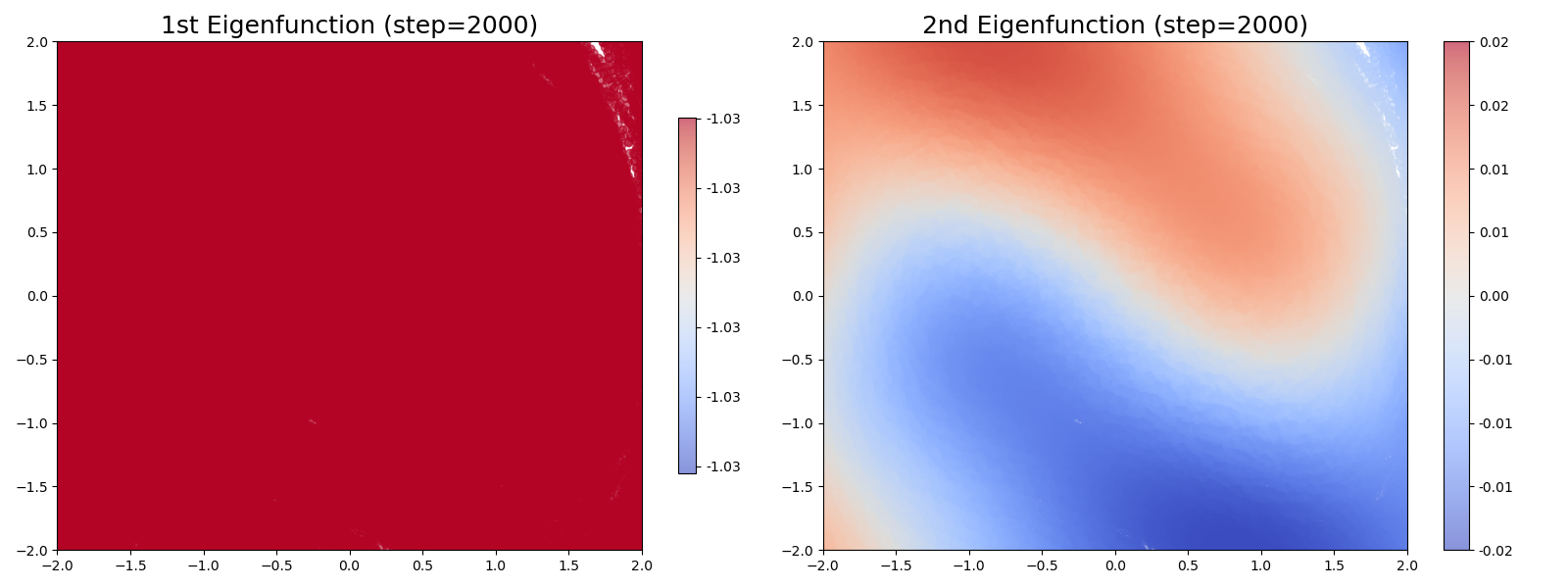}
        \caption{Step 2000}
        \label{fig:duffing_dqn_step_2000}
    \end{subfigure}

    \caption{Stochastic Duffing oscillator: the first two Koopman eigenfunctions learned by DQN-SDMD.}
    \label{fig:duffing_dqn_eigenfunctions}

\end{figure}

with parameters $\delta=0.5$, $\alpha=-1$, $\beta=1$, and $\sigma=0.15$. The deterministic drift has two stable equilibria at $(\pm1,0)$ separated by the saddle point $(0,0)$, as shown in Figure~\ref{fig:duffing_phase_portrait}. Damping drives trajectories toward one of the two wells, while stochastic forcing can induce rare basin crossings. Figure~\ref{fig:duffing_dqn_eigenfunctions} shows the learned Koopman eigenfunctions. The leading eigenfunction approaches the invariant mode. The
second eigenfunction captures the slowest non-trivial process, namely noise-induced switching between the two attractors. It is nearly constant inside each basin and changes sign across the separatrix.

\subsection{Stochastic FitzHugh--Nagumo Model}
\begin{wrapfigure}{r}{0.35\linewidth}
    \centering
    \includegraphics[width=0.95\linewidth]{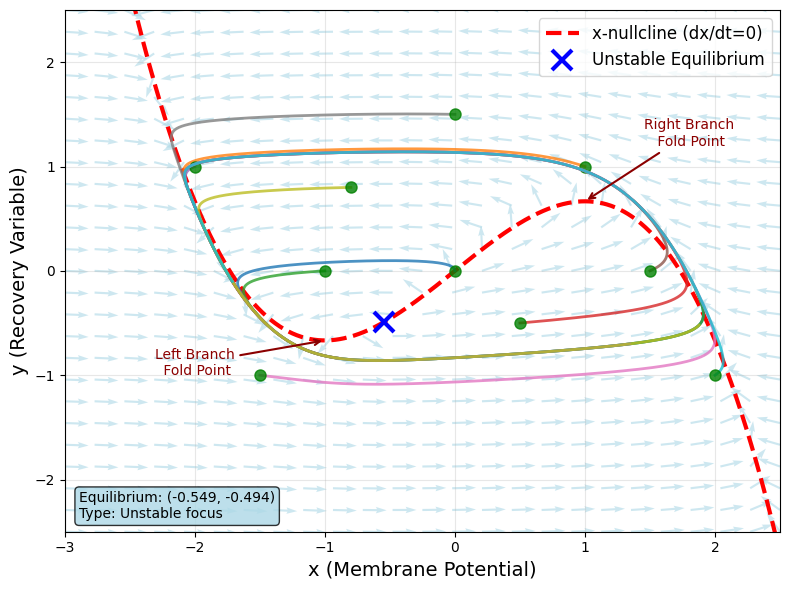}
    \caption{Stochastic FitzHugh--Nagumo model: phase portrait.}
    \label{fig:fhn_phase_portrait}
    \vspace{-2em}
\end{wrapfigure}

We next analyze the stochastic FitzHugh--Nagumo (FHN) model~\cite{fitzhugh1961impulses,nagumo1962active}, a prototypical fast--slow model of excitable neuronal dynamics. The fast variable $x$ represents the membrane potential and the slow variable $y$ represents the recovery current:
\begin{equation*}
    \begin{cases}
    \mathrm{d}x(t) = \bigl(x(t) - \tfrac{1}{3}x(t)^3 - y(t)\bigr)\,\mathrm{d}t + \sigma_1\,\mathrm{d}W_1(t), \\
    \mathrm{d}y(t) = \varepsilon \bigl(x(t) + a_1 - a_2 y(t)\bigr)\,\mathrm{d}t + \sigma_2\,\mathrm{d}W_2(t),
    \end{cases}
\end{equation*}
with parameters $\varepsilon=0.01$, $a_1=0.5$, $a_2=0.1$, $\sigma_1=10^{-3}$, and $\sigma_2=10^{-5}$. The small parameter $\varepsilon$ separates the fast $x$-dynamics from the slow $y$-dynamics and produces relaxation-type oscillations, as shown in Figure~\ref{fig:fhn_phase_portrait}. Trajectories drift slowly along the branches of the cubic $x$-nullcline and then undergo rapid transitions near the fold points $x\approx\pm1$.




Figure~\ref{fig:fhn_efuns} shows that the second non-trivial eigenfunction identifies the two long-residence parts of the phase space, near the top-left and bottom-right branches of the oscillation. These regions are mainly governed by the slow recovery variable, while the middle region corresponds to fast transitions between branches.

The unstable fixed point near $(-0.549,-0.494)$ appears as a low-density region, and the grey strip near the nullcline separates the two phases. For phase-reduced noisy limit-cycle dynamics, generator eigenvalues often have the form $\lambda_n=-n^2D_\phi+i n\omega$ with phase-diffusion coefficient $D_\phi>0$~\cite{kato2021asymptotic,takata2023definition}.

\begin{figure}[!t]
    \centering
    \begin{subfigure}{0.65\linewidth}
        \centering
        \includegraphics[width=\linewidth]{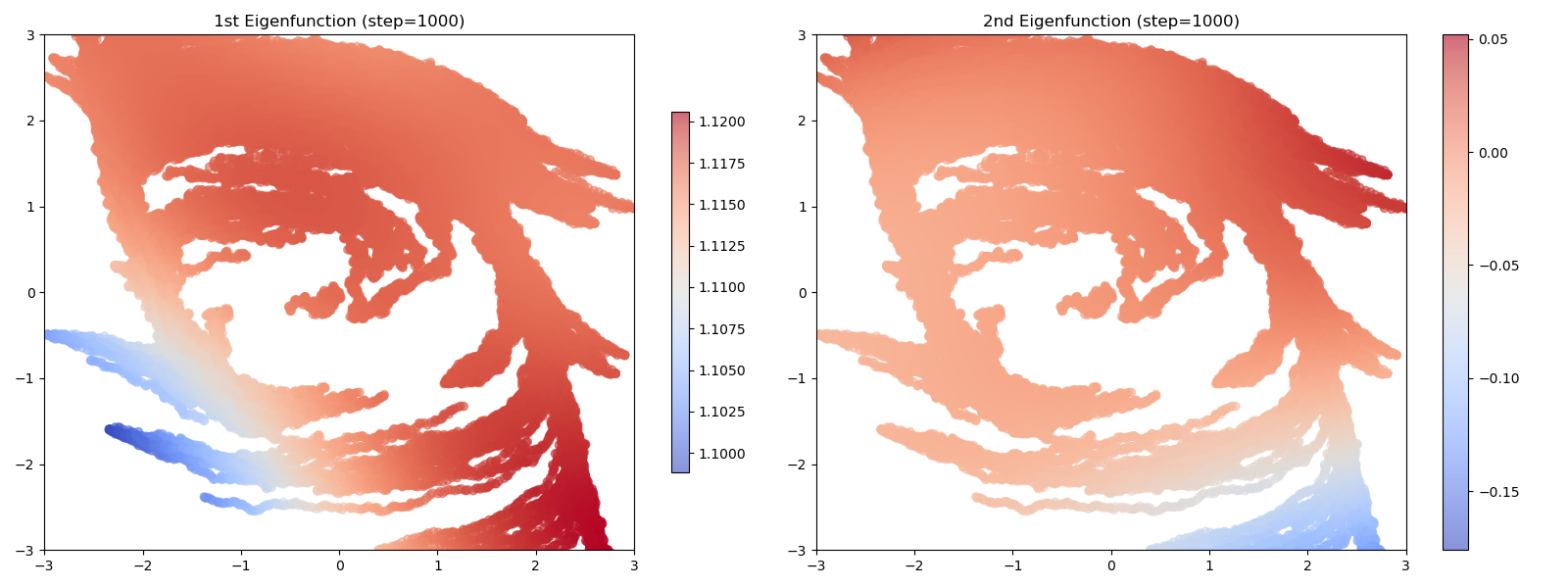}
        \caption{1st \& 2nd eigenfunctions at step 1000}
        \label{fig:eigen_1000}
    \end{subfigure}
    \medskip
    \begin{subfigure}{0.65\linewidth}
        \centering
        \includegraphics[width=\linewidth]{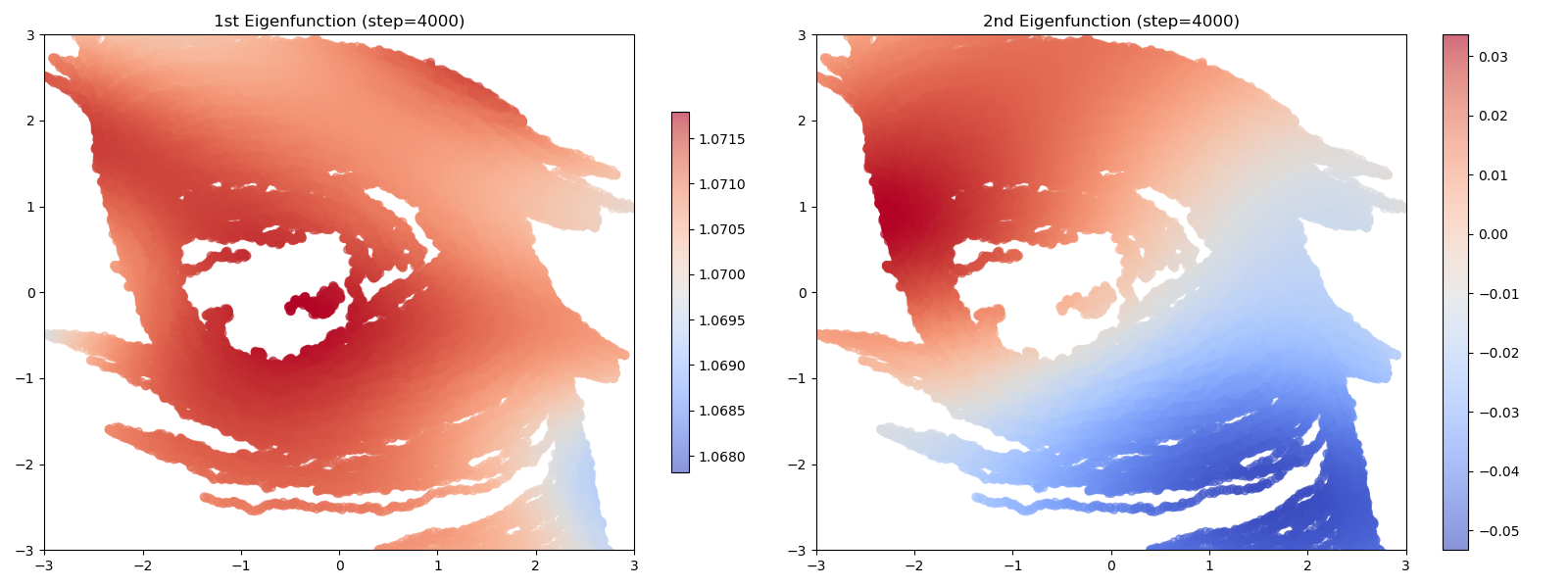}
        \caption{1st \& 2nd eigenfunctions at step 4000}
        \label{fig:eigen_4000}
    \end{subfigure}

    \medskip

    \begin{subfigure}{0.45\linewidth}
        \centering
        \includegraphics[width=\linewidth]{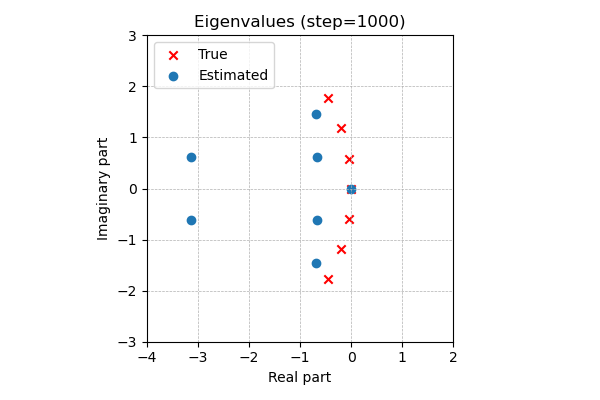}
        \caption{Step 1000}
        \label{fig:evalue_1000}
    \end{subfigure}
    \hspace{-1.5cm}
    \begin{subfigure}{0.45\linewidth}
        \centering
        \includegraphics[width=\linewidth]{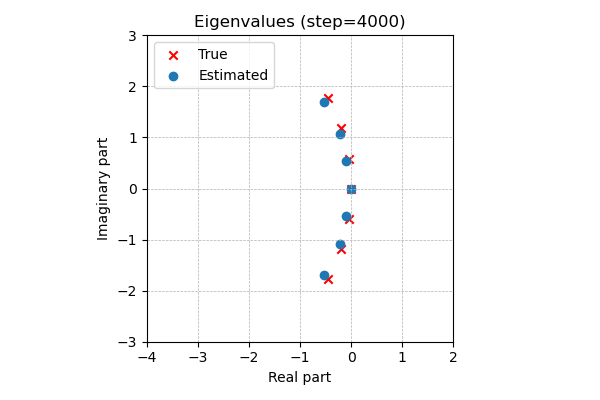}
        \caption{Step 4000}
        \label{fig:evalue_4000}
    \end{subfigure}
    \caption{First two eigenfunctions and generator eigenvalues in the larger-noise FitzHugh--Nagumo regime with $\sigma_1=\sigma_2=0.5$.}
    \label{fig:extra_fhn_images}
\end{figure}

Figure~\ref{fig:fhn_evalues} shows the corresponding evolution of the leading estimated eigenvalues. We also consider a larger-noise regime with $\sigma_1=\sigma_2=0.5$. In this case the limit-cycle structure is less clear, but the learned eigenfunctions still show the unstable fixed point and the separation between the two dominant phases; see Figure~\ref{fig:extra_fhn_images}.

\newcommand{\fhnsubfigwidth}{0.32}
\begin{figure*}[!t]
    \centering
    \begin{subfigure}{\fhnsubfigwidth\linewidth}
        \centering
        \includegraphics[width=\linewidth]{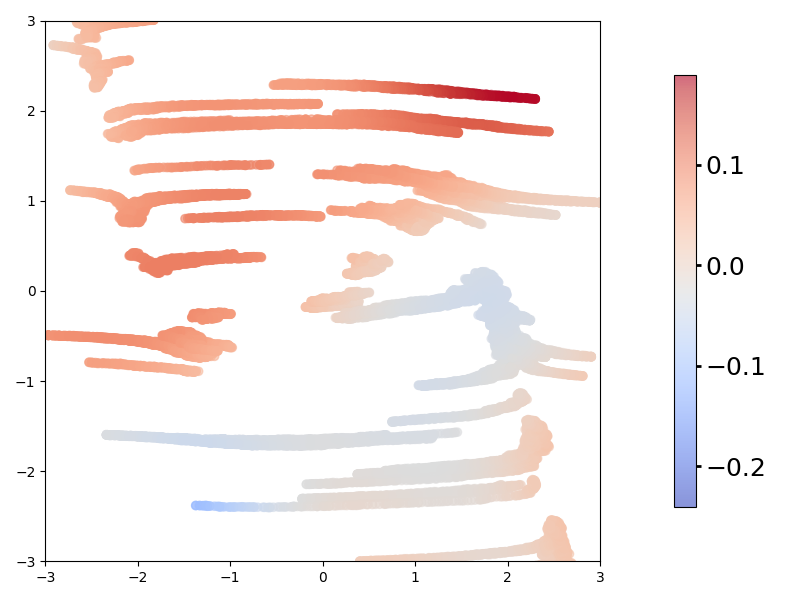}
        \caption{Step 500}
        \label{fig:fhn_efuns500}
    \end{subfigure}
    \hfill
    \begin{subfigure}{\fhnsubfigwidth\linewidth}
        \centering
        \includegraphics[width=\linewidth]{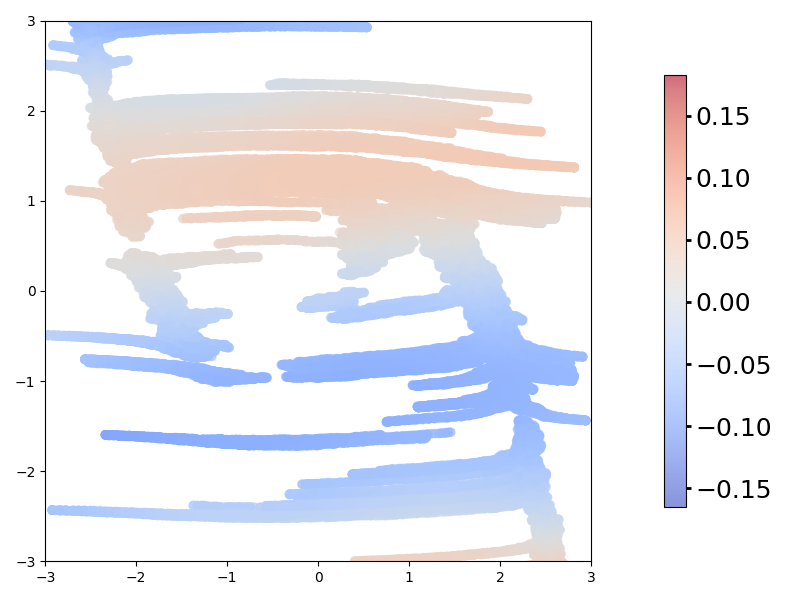}
        \caption{Step 2000}
        \label{fig:fhn_efuns2000}
    \end{subfigure}
    \hfill
    \begin{subfigure}{\fhnsubfigwidth\linewidth}
        \centering
        \includegraphics[width=\linewidth]{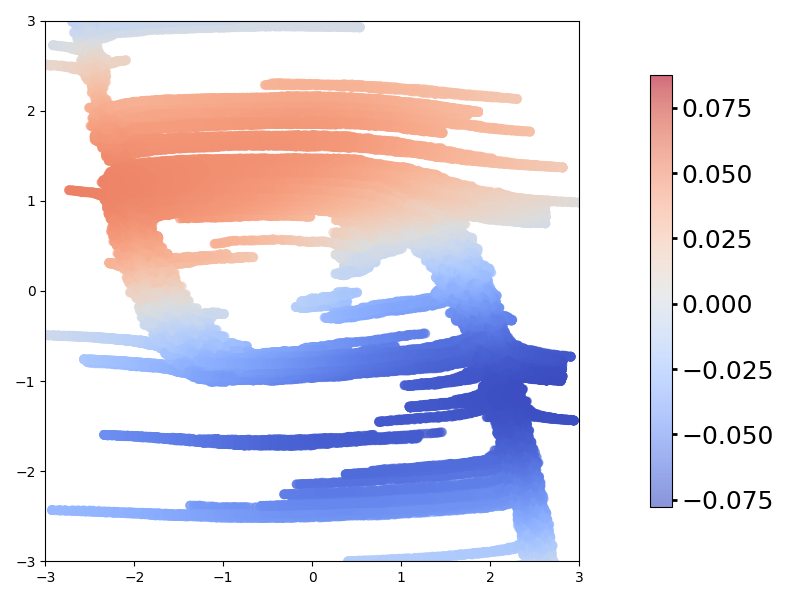}
        \caption{Step 8000}
        \label{fig:fhn_efuns8000}
    \end{subfigure}
    \caption{Second approximated Koopman eigenfunction of the FitzHugh--Nagumo model at different training steps.}
    \label{fig:fhn_efuns}
\end{figure*}

\begin{figure*}[!t]
    \centering
    \begin{subfigure}{\fhnsubfigwidth\linewidth}
        \centering
        \includegraphics[width=\linewidth]{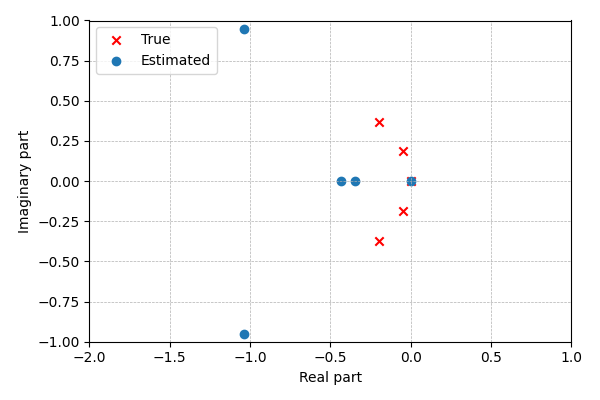}
        \caption{Step 500}
        \label{fig:fhn_step500}
    \end{subfigure}
    \hfill
    \begin{subfigure}{\fhnsubfigwidth\linewidth}
        \centering
        \includegraphics[width=\linewidth]{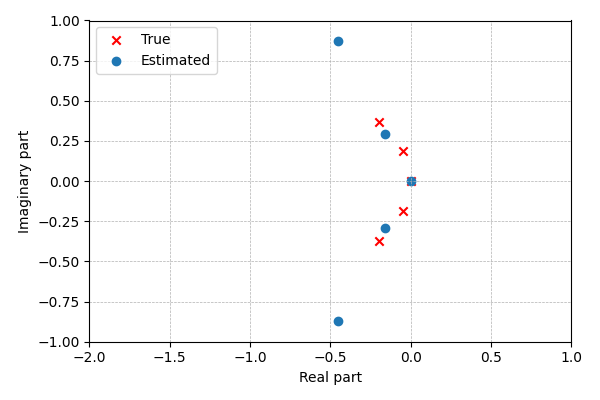}
        \caption{Step 2000}
        \label{fig:fhn_step2000}
    \end{subfigure}
    \hfill
    \begin{subfigure}{\fhnsubfigwidth\linewidth}
        \centering
        \includegraphics[width=\linewidth]{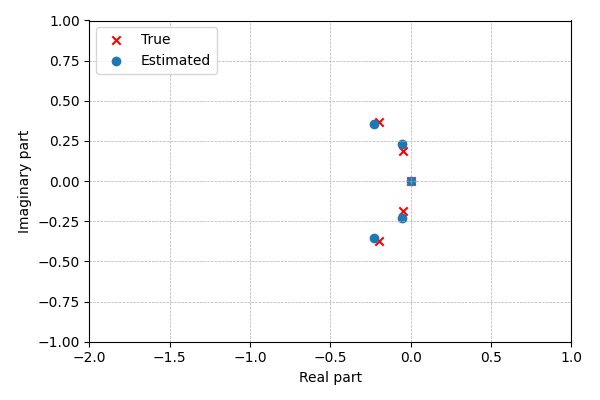}
        \caption{Step 8000}
        \label{fig:fhn_step8000}
    \end{subfigure}
    \caption{Approximated Koopman eigenvalues of the FitzHugh--Nagumo model at different training steps.}
    \label{fig:fhn_evalues}
\end{figure*}

\section{Convergence and Error Analysis}
This section gives error-propagation bounds for Bandit-SDMD, DQN-SDMD, and PPO-SDMD. The results should be read as idealized bounds: they do not prove finite-sample convergence of neural-network DQN or PPO training. Instead, they show how SDMD operator error enters standard bandit, approximate value-iteration, and approximate policy-iteration analyses.

We work with an infinite-horizon discounted MDP $(\mathsf{S},\mathsf{A},P,R,\gamma)$, where $0<\gamma<1$. Let $J(\pi)$ denote the expected discounted return of policy $\pi$ from a fixed initial distribution. Koopman operators act on bounded observables by conditional expectation. All operator norms below are induced norms on a bounded observable class containing the Bellman targets used by the algorithms.

\begin{assumption}[Uniform SDMD operator error]\label{ass:sdmd_error}
Let $\mathcal{K}$ denote an action-conditioned or policy-conditioned Koopman operator, and let $\widehat{\mathcal{K}}$ be its SDMD approximation on the chosen finite-dimensional feature space. We assume that
\begin{equation*}
    \|\mathcal{K}\|\le1,
    \qquad
    \|\widehat{\mathcal{K}}\|\le1,
    \qquad
    \|\mathcal{K}-\widehat{\mathcal{K}}\|\le\varepsilon_{\mathrm{sdmd}}.
\end{equation*}
Here $\varepsilon_{\mathrm{sdmd}}$ is the finite-data, fixed-dictionary, fixed-sampling-time SDMD operator error. The non-expansiveness assumption keeps the constants simple. If the finite-dimensional operators are only uniformly bounded, the same arguments hold with the corresponding bound constants.
\end{assumption}

\subsection{Regret Analysis of Bandit-SDMD}
For Bandit-SDMD, each arm $a\in\mathsf{A}$ is an initialization region with mean reward $R_a$. Let $R^*=\max_a R_a$ and define $\Delta_a=R^*-R_a$. The cumulative regret is $\mathcal{R}_T=\sum_{t=1}^T\Delta_{a_t}$.

\begin{assumption}[Reward-estimation precision]\label{ass:sdmd_precision}
The action values used for greedy selection satisfy
\begin{equation*}
|Q_t(a)-R_a|\le \varepsilon_{\mathrm{rew}},
\qquad \forall a\in\mathsf{A},\ \forall t.
\end{equation*}
Let $\Delta_{\min}=\min\{\Delta_a:\Delta_a>0\}$, with the convention that the condition is trivial if all arms are optimal. We assume
\begin{equation*}
    \varepsilon_{\mathrm{rew}}<\Delta_{\min}/2.
\end{equation*}
This assumption may be interpreted deterministically or on a high-probability event after the empirical rewards have been estimated accurately enough. If the reward is locally Lipschitz as a function of the estimated finite-dimensional Koopman matrix and the eigenvalues used in the reward are simple and separated, then one may take $\varepsilon_{\mathrm{rew}}\le L_R\varepsilon_{\mathrm{sdmd}}$ for a suitable perturbation constant $L_R$.
\end{assumption}

\begin{theorem}[Regret bound for Bandit-SDMD]\label{thm:sublinear_regret}
Consider an $\varepsilon$-greedy Bandit-SDMD algorithm with exploration schedule $\varepsilon_t=\min\{1,c/(|\mathsf{A}|t)\}$ for some $c>0$. On the event in Assumption~\ref{ass:sdmd_precision},
\begin{equation*}
\mathbb{E}[\mathcal{R}_T]
\le
\left(\frac{1}{|\mathsf{A}|}\sum_{a\in\mathsf{A}}\Delta_a\right)
\sum_{t=1}^T\varepsilon_t
=O(\log T),
\end{equation*}
where the expectation is over the exploration randomness. Hence $\mathbb{E}[\mathcal{R}_T]/T\to0$. If the precision event has probability at least $1-\delta$, then the unconditional bound gains the additional term $\delta T\Delta_{\max}$, where $\Delta_{\max}=\max_a\Delta_a$.
\end{theorem}
\begin{proof} See Appendix~\ref{thmpf:sublinear_regret}. \end{proof}

\begin{remark}
If the separation condition fails, a near-optimal arm may be indistinguishable from an optimal arm under the SDMD-induced reward error. Then a greedy step can repeatedly select the wrong arm, and linear regret is possible in the worst case; see Appendix~\ref{app:failure_case_bandit}.
\end{remark}

\subsection{Performance Bound for DQN-SDMD}
We analyze DQN-SDMD as an approximate value-iteration method. This separates the SDMD error from optimization, replay, and minibatch errors.

\begin{assumption}[Linear state--action representation]\label{ass:linear_q}
For every stationary policy $\pi$, the action-value function lies in the span of a bounded feature map $\mathbf{g}(s,a)$, i.e., $Q^\pi(s,a)=w_\pi^T\mathbf{g}(s,a)$. In addition, each function $f_k(s)=\max_{a'}Q_k(s,a')$ used in the Bellman optimality update belongs to the observable class in Assumption~\ref{ass:sdmd_error}.
\end{assumption}

Under this assumption, Lemma~\ref{lemma:q_target_error} shows that replacing the exact transition expectation by the SDMD approximation changes one Bellman target by at most $\gamma Q_{\max}\varepsilon_{\mathrm{sdmd}}$.

\begin{theorem}[Performance bound for DQN-SDMD]\label{thm:convergence_dqn}
Assume $\|Q_k\|_\infty\le Q_{\max}$ and suppose the DQN-SDMD update can be written as
\begin{equation*}
Q_{k+1}=\widehat{\mathcal{T}}Q_k+e_k,
\qquad
\limsup_{k\to\infty}\|e_k\|_\infty\le \varepsilon_{\mathrm{opt}},
\end{equation*}
where $\mathcal{T}$ is the exact Bellman optimality operator, $\widehat{\mathcal{T}}$ is the SDMD-induced approximate Bellman operator, and $e_k$ contains optimization and sampling error. If $\pi_k$ is greedy with respect to $Q_k$, then
\begin{equation}\label{eq:j_error_bound}
\limsup_{k\to\infty}\bigl(J(\pi^*)-J(\pi_k)\bigr)
\le
\frac{2\bigl(\gamma Q_{\max}\varepsilon_{\mathrm{sdmd}}+\varepsilon_{\mathrm{opt}}\bigr)}{(1-\gamma)^2}.
\end{equation}
\end{theorem}
\begin{proof} See Appendix~\ref{thmpf:convergence_dqn}. \end{proof}

\subsection{Performance Bound for PPO-SDMD}
We analyze PPO-SDMD as inexact approximate policy iteration. SDMD gives an approximate policy-evaluation step, and the PPO update is treated as an approximate policy-improvement step.

\begin{assumption}[Linear value and reward representation]\label{ass:linear_v}
For every stationary policy $\pi$, there are vectors $w_\pi,r_\pi\in\mathbb{R}^N$ such that
\[
    V^\pi(s)=w_\pi^T\mathbf{g}(s),
    \qquad
    \mathbb{E}_{a\sim\pi(\cdot\mid s)}[R(s,a)]=r_\pi^T\mathbf{g}(s).
\]
Assume $\|w_\pi\|\le W_{\max}$ and $C_g=\sup_{s\in\mathsf{S}}\|\mathbf{g}(s)\|<\infty$.
\end{assumption}

Lemma~\ref{lemma:value_error} gives
\[
    \|\widehat V^\pi-V^\pi\|_\infty
    \le
    \frac{\gamma C_g W_{\max}}{1-\gamma}\varepsilon_{\mathrm{sdmd}}.
\]
Combining this value-estimation bound with a standard approximate policy-iteration bound gives the following result.

\begin{theorem}[Performance bound for PPO-SDMD]\label{thm:convergence_ppo}
Suppose PPO-SDMD is modeled as an inexact approximate policy-iteration scheme with policy-improvement error bounded by $\varepsilon_{\mathrm{ppo}}$. Under Assumptions~\ref{ass:sdmd_error} and~\ref{ass:linear_v},
\begin{equation}\label{eq:ppo_error_bound}
\limsup_{k\to\infty}\bigl(J(\pi^*)-J(\pi_k)\bigr)
\le
\frac{2\gamma}{(1-\gamma)^2}
\left(
\frac{\gamma C_g W_{\max}}{1-\gamma}\varepsilon_{\mathrm{sdmd}}
+\varepsilon_{\mathrm{ppo}}
\right).
\end{equation}
Thus the SDMD contribution is $O(\varepsilon_{\mathrm{sdmd}}/(1-\gamma)^3)$, up to constants and the separate PPO improvement error.
\end{theorem}
\begin{proof} See Appendix~\ref{thmpf:convergence_ppo}. \end{proof}

\begin{remark}
The bounds separate the Koopman-estimation error from the RL optimization error. Better trajectory selection is intended to reduce $\varepsilon_{\mathrm{sdmd}}$, while network architecture, batch size, and training quality affect $\varepsilon_{\mathrm{opt}}$ or $\varepsilon_{\mathrm{ppo}}$.
\end{remark}

\section{Conclusion}
We introduced Reinforced SDMD, a method for learning where to collect trajectory data for stochastic Koopman spectral estimation. The method treats trajectory initialization as a sequential design problem: an RL agent selects initialization regions, SDMD estimates the Koopman semigroup, and a spectral-consistency reward guides the next sampling decisions.

The double-well, Duffing, and FitzHugh--Nagumo experiments show that bandit-, DQN-, and PPO-based sampling policies can concentrate data collection in dynamically relevant regions, including stable wells, basins of attraction, transition regions, and fast--slow phase partitions. The error analysis explains how the final learning performance depends on SDMD operator error together with the approximation and optimization errors of the chosen RL method.

Several directions remain open. A continuous-action version could replace the grid by a policy over initial conditions or over initialization distributions. Dictionary learning and policy learning could also be coupled more tightly, so that the observables and the sampling policy improve together. Finally, applications to fluid dynamics, molecular systems, and computational neuroscience would be valuable, especially in settings where short trajectories can be restarted from selected configurations. For fully passive datasets, the method is better viewed as adaptive subsampling or reweighting rather than closed-loop trajectory generation.

\section*{Acknowledgments}
The authors are grateful to Kaidi Shao and Rory Bufacchi for valuable early discussions, and to Parham Mohammad Panahi for helpful discussions on reinforcement-learning background. Y.X. and I.I. acknowledge support from JST CREST Grant Number JPMJCR24Q1, including the AIP Challenge Program. Z.S. was partially supported by NSERC RGPIN-2024-04938.

\bibliographystyle{plainnat}
\bibliography{references}

\newpage
\appendix
\titleformat{\section}[hang]{\normalfont\Large\bfseries}{}{0pt}{}

\section{Additional Algorithmic Details}


\subsection{Detailed Review of RL Algorithms}\label{app:rl_review}
This appendix records the RL components used in Reinforced SDMD. The notation follows Section~3.

\paragraph{Multi-armed bandit with $\varepsilon$-greedy policy.}
In the bandit setting~\cite{slivkins2019introduction}, each action $a$ selects one spatial region for trajectory initialization. After receiving reward $R_t$ and incrementing $N(a_t)$, the empirical action value is updated by
\begin{equation*}
Q(a_t)\leftarrow Q(a_t)+\frac{1}{N(a_t)}\bigl(R_t-Q(a_t)\bigr),
\end{equation*}
where $N(a_t)$ is the number of times action $a_t$ has been selected. The $\varepsilon$-greedy rule selects the current best action with probability $1-\varepsilon$ and a uniformly random action otherwise.

\paragraph{Deep Q-Network (DQN).}
DQN~\cite{mnih2013playing} uses a state $S_t$, which we take to be a recent history of initialization regions. A neural network approximates $Q(S_t,a_t;\theta)$ and is trained through the temporal-difference error
\begin{equation*}
\delta_t=R_t+\gamma\max_{a'}Q(S_{t+1},a';\theta^-)-Q(S_t,a_t;\theta),
\end{equation*}
where $\theta^-$ denotes the target-network parameters. Experience replay reduces temporal correlation in minibatches, and target-network updates stabilize the bootstrapped target.

\paragraph{Proximal Policy Optimization (PPO).}
PPO~\cite{schulman2017proximal} learns a stochastic policy $\pi_\theta(a\mid s)$ together with a critic $V(s)$. Its clipped surrogate objective is
\begin{equation}\label{eq:ppo_clip}
L_{\mathrm{CLIP}}(\theta)
=
\mathbb{E}_t\!\left[
\min\!\Bigl(r_t(\theta)\widehat A_t,
\operatorname{clip}\bigl(r_t(\theta),1-\varepsilon_{\mathrm{clip}},1+\varepsilon_{\mathrm{clip}}\bigr)\widehat A_t\Bigr)
\right],
\end{equation}
where $r_t(\theta)=\pi_\theta(a_t\mid s_t)/\pi_{\theta_{\old}}(a_t\mid s_t)$ and $\widehat A_t$ is an advantage estimate. We use generalized advantage estimation (GAE),
\begin{equation*}
\widehat A_t=\sum_{l=0}^{T-t-1}(\gamma\lambda)^l\delta_{t+l},
\qquad
\delta_t=R_t+\gamma V(s_{t+1})-V(s_t).
\end{equation*}
The clipping term limits the size of the policy update, which is useful here because the spectral reward changes as the Koopman estimate is updated.

\subsection{Failure of the Reward-Separation Condition}\label{app:failure_case_bandit}
When Assumption~\ref{ass:sdmd_precision} is violated, the estimated action values may be unable to distinguish an optimal arm from a near-optimal suboptimal arm. A mistaken greedy choice of a suboptimal arm $a$ is possible only when $\Delta_a\le 2\varepsilon_{\mathrm{rew}}$. Therefore, if $2\varepsilon_{\mathrm{rew}}>\Delta_{\min}$, with $\Delta_{\min}=\min\{\Delta_a:\Delta_a>0\}$, the set of potentially confusable suboptimal arms may be nonempty.

A simple worst-case regret bound is
\begin{equation*}
\mathbb{E}[\mathcal{R}_T]
\le
\sum_{t=1}^T\varepsilon_t\left(\frac{1}{|\mathsf{A}|}\sum_{a\in\mathsf{A}}\Delta_a\right)
+
\sum_{t=1}^T(1-\varepsilon_t)
\left(\sum_{a:\,\Delta_a\le 2\varepsilon_{\mathrm{rew}}}\Delta_a\right).
\end{equation*}
The first term is the exploration cost. The second term upper-bounds the exploitation loss when the estimator bias persistently favors a confusable suboptimal arm. Let
\[
C_{\mathrm{err}}\coloneqq\sum_{a:\,\Delta_a\le2\varepsilon_{\mathrm{rew}}}\Delta_a.
\]
For $\varepsilon_t=c/(|\mathsf{A}|t)$, the exploration term is $O(\log T)$ whereas the exploitation term can scale as $C_{\mathrm{err}}(T-O(\log T))$. Hence the regret can be $\Theta(T)$ in the worst case.

\subsection{DQN Implementation Details}\label{app:dqn}
The DQN experiments use the following standard stabilization mechanisms.

\paragraph{Huber loss.}
Instead of mean squared error, we use the Huber loss~\cite{huber1992robust},
\begin{equation*}
L_\delta(y,\widehat y)=
\begin{cases}
\frac12(y-\widehat y)^2, & |y-\widehat y|\le\delta,\\[4pt]
\delta\left(|y-\widehat y|-\frac12\delta\right), & \text{otherwise.}
\end{cases}
\end{equation*}
This loss is quadratic for small TD errors and linear for large TD errors, making it more robust to outliers during early exploration.

\paragraph{Experience replay.}
Transitions are stored in a replay buffer and minibatches are sampled randomly for training. This reduces correlation between consecutive updates and improves stability.

\paragraph{Soft target-network updates.}
The target network is updated by $\theta^-\leftarrow\tau\theta+(1-\tau)\theta^-$, with $\tau=0.05$ in the reported implementation.

\paragraph{Gradient clipping.}
Gradients are clipped to avoid exploding updates. In the implementation this is performed by \texttt{torch.nn.utils.clip\_grad\_value\_(policy\_net.parameters(), 100)}.

\subsection{Clipped Surrogate Objective in PPO}
The PPO clipped surrogate objective controls the change from the old policy to the new policy through the likelihood ratio
\[
r_t(\theta)=\frac{\pi_\theta(a_t\mid S_t)}{\pi_{\theta_{\old}}(a_t\mid S_t)}.
\]
The unclipped policy-gradient objective is $L_{\mathrm{PG}}(\theta)=\mathbb{E}_t[r_t(\theta)\widehat A_t]$, which may produce unstable updates when $r_t$ is too large or too small. PPO replaces it by the clipped objective in Eq.~\eqref{eq:ppo_clip}. When $\widehat A_t>0$, clipping prevents the new policy from increasing the probability of the action by more than a factor $1+\varepsilon_{\mathrm{clip}}$; when $\widehat A_t<0$, clipping prevents the probability from being reduced below the factor $1-\varepsilon_{\mathrm{clip}}$.

\subsection{Sufficient Conditions for Assumptions~\ref{ass:linear_q} and~\ref{ass:linear_v}}\label{app:w_derivation_detail}
This section states sufficient conditions for the linear value-function assumptions used in Section~5.

\paragraph{Value-function case.}
For a fixed policy $\pi$, the Bellman equation is
\begin{equation*}
V^\pi(s)=\mathbb{E}_{a\sim\pi(\cdot\mid s),\,s'\sim P(\cdot\mid s,a)}
\bigl[R(s,a)+\gamma V^\pi(s')\bigr].
\end{equation*}
Assume $V^\pi(s)=w_\pi^T\mathbf g(s)$ and $\mathbb{E}_{a\sim\pi(\cdot\mid s)}[R(s,a)]=r_\pi^T\mathbf g(s)$. If the feature span is invariant under $\mathcal{K}_\pi$, so that $(\mathcal{K}_\pi\mathbf g)(s)=M_\pi^T\mathbf g(s)$, then
\[
w_\pi^T\mathbf g(s)=\bigl(r_\pi+\gamma M_\pi w_\pi\bigr)^T\mathbf g(s),
\]
and hence $(I-\gamma M_\pi)w_\pi=r_\pi$. Since $\|M_\pi\|\le1$, the inverse exists for $\gamma<1$ and
\[
(I-\gamma M_\pi)^{-1}=\sum_{n=0}^{\infty}(\gamma M_\pi)^n,
\qquad
\|(I-\gamma M_\pi)^{-1}\|\le\frac{1}{1-\gamma}.
\]

\paragraph{Q-function case.}
For a fixed policy $\pi$, the state--action process $(S_t,A_t)$ is Markov. Define the Koopman operator on state--action observables by
\[
(\mathcal{K}_\pi f)(s,a)
=\mathbb{E}\bigl[f(S_{t+1},A_{t+1})\mid S_t=s,A_t=a,
A_{t+1}\sim\pi(\cdot\mid S_{t+1})\bigr].
\]
Then $Q^\pi=R+\gamma\mathcal{K}_\pi Q^\pi$. If the reward is representable as $R(s,a)=r^T\mathbf g(s,a)$ and the feature span is invariant under $\mathcal{K}_\pi$, then $Q^\pi(s,a)=w_\pi^T\mathbf g(s,a)$ for a vector $w_\pi$ solving the corresponding finite-dimensional linear system.

\section{Proofs}

\subsection{Proof of Theorem~\ref{thm:sublinear_regret}}\label{thmpf:sublinear_regret}
\begin{proof}
Work on the precision event in Assumption~\ref{ass:sdmd_precision}. At time $t$, regret can come from random exploration or from a wrong greedy choice. The greedy choice cannot be suboptimal on this event. Indeed, if a suboptimal arm $a$ were selected greedily, then $Q_t(a)\ge Q_t(a^*)$ for some optimal arm $a^*$. Hence
\[
R_a+\varepsilon_{\mathrm{rew}}
\ge Q_t(a)
\ge Q_t(a^*)
\ge R^*-\varepsilon_{\mathrm{rew}},
\]
so $\Delta_a=R^*-R_a\le2\varepsilon_{\mathrm{rew}}$. This contradicts $\varepsilon_{\mathrm{rew}}<\Delta_{\min}/2$ for any suboptimal arm.

Thus, on the precision event, exploitation adds no regret. Exploration chooses uniformly over $\mathsf{A}$, and therefore
\[
\mathbb{E}[\mathcal{R}_T\mid \mathcal{E}]
\le
\sum_{t=1}^T\varepsilon_t
\left(\frac{1}{|\mathsf{A}|}\sum_{a\in\mathsf{A}}\Delta_a\right).
\]
Since $\varepsilon_t=\min\{1,c/(|\mathsf{A}|t)\}$, the sum $\sum_{t=1}^T\varepsilon_t$ is $O(\log T)$. This proves the stated conditional bound. If $\mathbb{P}(\mathcal{E})\ge1-\delta$, then on $\mathcal{E}^c$ the regret is at most $T\Delta_{\max}$, which adds the term $\delta T\Delta_{\max}$.
\end{proof}

\subsection{Bellman Target Error}\label{lemmapf:q_target_error}
\begin{lemma}[Bellman target estimation error]\label{lemma:q_target_error}
Let $\mathcal{T}$ be the exact Bellman optimality operator and $\widehat{\mathcal{T}}$ be the operator obtained by replacing the transition expectation with the SDMD-estimated Koopman expectation. Under Assumptions~\ref{ass:sdmd_error} and~\ref{ass:linear_q}, if $\|Q_k\|_\infty\le Q_{\max}$, then
\begin{equation*}
\|\mathcal{T}Q_k-\widehat{\mathcal{T}}Q_k\|_\infty
\le
\gamma Q_{\max}\varepsilon_{\mathrm{sdmd}}.
\end{equation*}
\end{lemma}

\begin{proof}
Let $f_k(s)=\max_{a'}Q_k(s,a')$. Then $\|f_k\|_\infty\le\|Q_k\|_\infty\le Q_{\max}$. The only difference between $\mathcal{T}$ and $\widehat{\mathcal{T}}$ is the conditional expectation term, hence
\[
\|\mathcal{T}Q_k-\widehat{\mathcal{T}}Q_k\|_\infty
=\gamma\sup_{s,a}\left|((\mathcal{K}_a-\widehat{\mathcal{K}}_a)f_k)(s)\right|.
\]
By Assumption~\ref{ass:linear_q}, $f_k$ belongs to the observable class on which the SDMD error bound holds. Therefore
\[
\|\mathcal{T}Q_k-\widehat{\mathcal{T}}Q_k\|_\infty
\le
\gamma\sup_a\|\mathcal{K}_a-\widehat{\mathcal{K}}_a\|\,\|f_k\|_\infty
\le
\gamma Q_{\max}\varepsilon_{\mathrm{sdmd}}.
\]
\end{proof}

\subsection{Proof of Theorem~\ref{thm:convergence_dqn}}\label{thmpf:convergence_dqn}
\begin{proof}
Let $\delta_k=\|Q^*-Q_k\|_\infty$. Since $\mathcal{T}$ is a $\gamma$-contraction and $Q_{k+1}=\widehat{\mathcal{T}}Q_k+e_k$,
\[
\begin{aligned}
\delta_{k+1}
&=\|\mathcal{T}Q^*-\widehat{\mathcal{T}}Q_k-e_k\|_\infty\\
&\le \|\mathcal{T}Q^*-\mathcal{T}Q_k\|_\infty
+\|\mathcal{T}Q_k-\widehat{\mathcal{T}}Q_k\|_\infty
+\|e_k\|_\infty\\
&\le \gamma\delta_k+\gamma Q_{\max}\varepsilon_{\mathrm{sdmd}}+\|e_k\|_\infty.
\end{aligned}
\]
Taking the limsup and using $\limsup_k\|e_k\|_\infty\le\varepsilon_{\mathrm{opt}}$ gives
\[
\limsup_{k\to\infty}\|Q^*-Q_k\|_\infty
\le
\frac{\gamma Q_{\max}\varepsilon_{\mathrm{sdmd}}+\varepsilon_{\mathrm{opt}}}{1-\gamma}.
\]
For a policy $\pi_k$ greedy with respect to $Q_k$, the standard performance-loss inequality yields
\[
J(\pi^*)-J(\pi_k)
\le
\frac{2}{1-\gamma}\|Q^*-Q_k\|_\infty.
\]
Combining the two inequalities proves Eq.~\eqref{eq:j_error_bound}.
\end{proof}

\subsection{Value Estimation Error}\label{lemmapf:value_error}
\begin{lemma}[Value estimation error bound]\label{lemma:value_error}
Under Assumptions~\ref{ass:sdmd_error} and~\ref{ass:linear_v}, the SDMD policy-evaluation error satisfies
\begin{equation*}
\|\widehat V^\pi-V^\pi\|_\infty
\le
\frac{\gamma C_g W_{\max}}{1-\gamma}\varepsilon_{\mathrm{sdmd}}.
\end{equation*}
\end{lemma}

\begin{proof}
On the feature span, write
\[
 w_\pi=(I-\gamma M_\pi)^{-1}r_\pi,
 \qquad
 \widehat w_\pi=(I-\gamma\widehat M_\pi)^{-1}r_\pi.
\]
Using $A^{-1}-B^{-1}=A^{-1}(B-A)B^{-1}$ with $A=I-\gamma\widehat M_\pi$ and $B=I-\gamma M_\pi$ gives
\[
\widehat w_\pi-w_\pi
=
\gamma(I-\gamma\widehat M_\pi)^{-1}(\widehat M_\pi-M_\pi)w_\pi.
\]
Because $\|\widehat M_\pi\|\le1$, the Neumann-series bound gives $\|(I-\gamma\widehat M_\pi)^{-1}\|\le(1-\gamma)^{-1}$. Hence
\[
\|\widehat w_\pi-w_\pi\|
\le
\frac{\gamma}{1-\gamma}\varepsilon_{\mathrm{sdmd}}\|w_\pi\|
\le
\frac{\gamma W_{\max}}{1-\gamma}\varepsilon_{\mathrm{sdmd}}.
\]
Finally,
\[
\|\widehat V^\pi-V^\pi\|_\infty
=\sup_s|(\widehat w_\pi-w_\pi)^T\mathbf g(s)|
\le
C_g\|\widehat w_\pi-w_\pi\|,
\]
which proves the claim.
\end{proof}

\subsection{Proof of Theorem~\ref{thm:convergence_ppo}}\label{thmpf:convergence_ppo}
\begin{proof}
In the approximate policy-iteration abstraction, the asymptotic performance loss is bounded by
\[
\limsup_{k\to\infty}(J(\pi^*)-J(\pi_k))
\le
\frac{2\gamma}{(1-\gamma)^2}(\tau+\varepsilon_{\mathrm{ppo}}),
\]
where $\tau$ is the uniform policy-evaluation error and $\varepsilon_{\mathrm{ppo}}$ is the policy-improvement error. Lemma~\ref{lemma:value_error} gives
\[
\tau\le\frac{\gamma C_g W_{\max}}{1-\gamma}\varepsilon_{\mathrm{sdmd}}.
\]
Substituting this value of $\tau$ proves Eq.~\eqref{eq:ppo_error_bound}.
\end{proof}

\end{document}